\documentclass[12pt]{article}

\begin{document}

\newcommand{\llra}{\longleftrightarrow}
\newcommand{\ra}{\rightarrow}
\newcommand{\alg}{ {\cal A} }
\newcommand{\EJ}{\hat{J}_{sl(2),\lambda}}
\newcommand{\hb}{\hbar}
\newcommand{\rrr}{{\bf R}}
\newcommand{\ccd}{{\cal T}}
\newcommand{\ggg}{{\cal G}}
\newcommand{\aaa}{{\cal A}}
\newcommand{\ddd}{{\cal D}}
\newcommand{\zzz}{{\bf Z}}
\newcommand{\ci}{\epsilon}
\newcommand{\ccc}{{\bf C}}
\newcommand{\zsn}{{\cal Z}S_n}
\newcommand{\lra}{\longrightarrow}
\newcommand{\bq}{\begin{equation}}
\newcommand{\eq}{\end{equation}}
\newcommand{\rtb}{\begin{flushright}
\begin{picture}(0.3,0.3)(0,0)
\put(-0.3,0){\framebox(0.3,0.3)}
\end{picture}
\end{flushright}}

\newtheorem{theorem}{Theorem}
\newtheorem{defn}{Definition}[section]
\newtheorem{lemma}[defn]{Lemma}
\newtheorem{corollary}[defn]{Corollary}
\newtheorem{fact}[defn]{Fact}
\newtheorem{prop}[defn]{Proposition}
\newtheorem{conj}[defn]{Conjecture}

%
%
%
%
%
\def\CabCirc{
\thicklines
\qbezier(0.966, -0.259)(1.04, 0)(0.966, 0.259)
\qbezier(0.966, 0.259)(0.897, 0.518)(0.707, 0.707)
\qbezier(0.707, 0.707)(0.518, 0.897)(0.259, 0.966)
\qbezier(-0.259, 0.966)(-0.518, 0.897)(-0.707, 0.707)
\qbezier(-0.707, 0.707)(-0.897, 0.518)(-0.966, 0.259)
\qbezier(-0.966, 0.259)(-1.04, 0)(-0.966, -0.259)
\qbezier(0.773, -0.207)(0.832 , 0)(0.773 , 0.207)
\qbezier(0.773, 0.207)(0.718, 0.414)(0.566, 0.566)
\qbezier(0.566, 0.566)(0.414, 0.718)(0.207, 0.773)
\qbezier(-0.207, 0.773)(-0.414, 0.718)(-0.566, 0.566)
\qbezier(-0.566, 0.566)(-0.718, 0.414)(-0.773, 0.207)
\qbezier(-0.773, 0.207)(-0.832, 0)(-0.773, -0.207)
\put(0,0){\oval(2,2)[b]}
\put(0,0){\oval(1.6,1.6)[b]}
\qbezier(-0.259,0.966)(.207,0.773)(.207,0.773)
\qbezier(-.207,0.773)(-0.207,0.773)(0.259,0.966)
}
\def\CabA{
\thinlines
\qbezier(-0.4,-0.693)(-0.4,-0.693)(0.5,0.866)
\put(-0.4,-0.693){\circle*{0.15}}
\put(0.5,0.866){\circle*{0.15}}
}
\def\CabB{
\thinlines
\qbezier(-0.5,-0.866)(-0.5,-0.866)(0.4,0.693)
\put(-0.5,-0.866){\circle*{0.15}}
\put(0.4,0.693){\circle*{0.15}}
}
\def\CabC{
\thinlines
\qbezier(-0.693,0.4)(-0.693,0.4)(0.866,-0.5)
\put(-0.693,0.4){\circle*{0.15}}
\put(0.866,-0.5){\circle*{0.15}}
}
\def\CabD{
\thinlines
\qbezier(-0.866,0.5)(0.693,-0.4)(0.693,-0.4)
\put(-0.866,0.5){\circle*{0.15}}
\put(0.693,-0.4){\circle*{0.15}}
}
\def\CabE{
\thinlines
\qbezier(0.693,-0.4)(0.693,-0.4)(-0.693,0.4)
\put(0.693,-0.4){\circle*{0.15}}
\put(-0.693,0.4){\circle*{0.15}}
}
\def\CabF{
\thinlines
\qbezier(0.4,0.693)(0.4,0.693)(-0.4,-0.693)
\put(0.4,0.693){\circle*{0.15}}
\put(-0.4,-0.693){\circle*{0.15}}
}

\setlength{\unitlength}{20pt}
\def\DottedCircle{
\qbezier[4](0.966,-0.259)(1.04,0)(0.966,0.259)
\qbezier[4](0.966,0.259)(0.897,0.518)(0.707,0.707)
\qbezier[4](0.707,0.707)(0.518,0.897)(0.259,0.966)
\qbezier[4](0.259,0.966)(0,1.04)(-0.259,0.966)
\qbezier[4](-0.259,0.966)(-0.518,0.897)(-0.707,0.707)
\qbezier[4](-0.707,0.707)(-0.897,0.518)(-0.966,0.259)
\qbezier[4](-0.966,0.259)(-1.04,0)(-0.966,-0.259)
\qbezier[4](-0.966,-0.259)(-0.897,-0.518)(-0.707,-0.707)
\qbezier[4](-0.707,-0.707)(-0.518,-0.897)(-0.259,-0.966)
\qbezier[4](-0.259,-0.966)(0,-1.04)(0.259,-0.966)
\qbezier[4](0.259,-0.966)(0.518,-0.897)(0.707,-0.707)
\qbezier[4](0.707,-0.707)(0.897,-0.518)(0.966,-0.259)
}
\def\FullCircle{
\thicklines
\put(0,0){\circle{2}}
}
%
%
\def\Endpoint[#1]{
\ifcase#1
\put(1,0){\circle*{0.15}}
\or\put(0.866,0.5){\circle*{0.15}}
\or\put(0.5,0.866){\circle*{0.15}}
\or\put(0,1){\circle*{0.15}}
\or\put(-0.5,0.866){\circle*{0.15}}
\or\put(-0.866,0.5){\circle*{0.15}}
\or\put(-1,0){\circle*{0.15}}
\or\put(-0.866,-0.5){\circle*{0.15}}
\or\put(-0.5,-0.866){\circle*{0.15}}
\or\put(0,-1){\circle*{0.15}}
\or\put(0.5,-0.866){\circle*{0.15}}
\or\put(0.866,-0.5){\circle*{0.15}}
\fi}
%
%
\def\Arc[#1]{
\ifcase#1
\qbezier[25](0.966,-0.259)(1.04,0)(0.966,0.259)
\or
\qbezier[25](0.966,0.259)(0.897,0.518)(0.707,0.707)
\or
\qbezier[25](0.707,0.707)(0.518,0.897)(0.259,0.966)
\or
\qbezier[25](0.259,0.966)(0,1.04)(-0.259,0.966)
\or
\qbezier[25](-0.259,0.966)(-0.518,0.897)(-0.707,0.707)
\or
\qbezier[25](-0.707,0.707)(-0.897,0.518)(-0.966,0.259)
\or
\qbezier[25](-0.966,0.259)(-1.04,0)(-0.966,-0.259)
\or
\qbezier[25](-0.966,-0.259)(-0.897,-0.518)(-0.707,-0.707)
\or
\qbezier[25](-0.707,-0.707)(-0.518,-0.897)(-0.259,-0.966)
\or
\qbezier[25](-0.259,-0.966)(0,-1.04)(0.259,-0.966)
\or
\qbezier[25](0.259,-0.966)(0.518,-0.897)(0.707,-0.707)
\or
\qbezier[25](0.707,-0.707)(0.897,-0.518)(0.966,-0.259)
\fi}
%
%
\def\Chord[#1,#2]{
\thinlines
\ifnum#1>#2\Chord[#2,#1]
\else\ifnum#1<#2
\ifcase#1
\ifcase#2
\or\qbezier(1,0)(0.516,0.138)(0.866,0.5)
\or\qbezier(1,0)(0.45,0.26)(0.5,0.866)
\or\qbezier(1,0)(0.327,0.327)(0,1)
\or\qbezier(1,0)(0.179,0.311)(-0.5,0.866)
\or\qbezier(1,0)(0.0536,0.2)(-0.866,0.5)
\or\put(1, 0){\line(-2, 0){2}}
\or\qbezier(1,0)(0.0536,-0.2)(-0.866,-0.5)
\or\qbezier(1,0)(0.179,-0.311)(-0.5,-0.866)
\or\qbezier(1,0)(0.327,-0.327)(0,-1)
\or\qbezier(1,0)(0.45,-0.26)(0.5,-0.866)
\or\qbezier(1,0)(0.516,-0.138)(0.866,-0.5)
\fi
\or\ifcase#2\or
\or\qbezier(0.866,0.5)(0.378,0.378)(0.5,0.866)
\or\qbezier(0.866,0.5)(0.26,0.45)(0,1)
\or\qbezier(0.866,0.5)(0.12,0.446)(-0.5,0.866)
\or\qbezier(0.866,0.5)(0,0.359)(-0.866,0.5)
\or\qbezier(0.866,0.5)(-0.0536,0.2)(-1,0)
\or\put(0.866, 0.5){\line(-5, -3){1.73}}
\or\qbezier(0.866,0.5)(0.146,-0.146)(-0.5,-0.866)
\or\qbezier(0.866,0.5)(0.311,-0.179)(0,-1)
\or\qbezier(0.866,0.5)(0.446,-0.12)(0.5,-0.866)
\or\qbezier(0.866,0.5)(0.52,0)(0.866,-0.5)
\fi
\or\ifcase#2\or\or
\or\qbezier(0.5,0.866)(0.138,0.516)(0,1)
\or\qbezier(0.5,0.866)(0,0.52)(-0.5,0.866)
\or\qbezier(0.5,0.866)(-0.12,0.446)(-0.866,0.5)
\or\qbezier(0.5,0.866)(-0.179,0.311)(-1,0)
\or\qbezier(0.5,0.866)(-0.146,0.146)(-0.866,-0.5)
\or\put(0.5, 0.866){\line(-3, -5){1}}
\or\qbezier(0.5,0.866)(0.2,-0.0536)(0,-1)
\or\qbezier(0.5,0.866)(0.359,0)(0.5,-0.866)
\or\qbezier(0.5,0.866)(0.446,0.12)(0.866,-0.5)
\fi
\or\ifcase#2\or\or\or
\or\qbezier(0,1.)(-0.138,0.516)(-0.5,0.866)
\or\qbezier(0,1.)(-0.26,0.45)(-0.866,0.5)
\or\qbezier(0,1.)(-0.327,0.327)(-1,0)
\or\qbezier(0,1.)(-0.311,0.179)(-0.866,-0.5)
\or\qbezier(0,1.)(-0.2,0.0536)(-0.5,-0.866)
\or\put(0, 1){\line(0, -2){2}}
\or\qbezier(0,1.)(0.2,0.0536)(0.5,-0.866)
\or\qbezier(0,1.)(0.311,0.179)(0.866,-0.5)
\fi
\or\ifcase#2\or\or\or\or
\or\qbezier(-0.5,0.866)(-0.378,0.378)(-0.866,0.5)
\or\qbezier(-0.5,0.866)(-0.45,0.26)(-1,0)
\or\qbezier(-0.5,0.866)(-0.446,0.12)(-0.866,-0.5)
\or\qbezier(-0.5,0.866)(-0.359,0)(-0.5,-0.866)
\or\qbezier(-0.5,0.866)(-0.2,-0.0536)(0,-1)
\or\put(-0.5, 0.866){\line(3, -5){1}}
\or\qbezier(-0.5,0.866)(0.146,0.146)(0.866,-0.5)
\fi
\or\ifcase#2\or\or\or\or\or
\or\qbezier(-0.866,0.5)(-0.516,0.138)(-1,0)
\or\qbezier(-0.866,0.5)(-0.52,0)(-0.866,-0.5)
\or\qbezier(-0.866,0.5)(-0.446,-0.12)(-0.5,-0.866)
\or\qbezier(-0.866,0.5)(-0.311,-0.179)(0,-1)
\or\qbezier(-0.866,0.5)(-0.146,-0.146)(0.5,-0.866)
\or\put(-0.866, 0.5){\line(5, -3){1.73}}
\fi
\or\ifcase#2\or\or\or\or\or\or
\or\qbezier(-1,0)(-0.516,-0.138)(-0.866,-0.5)
\or\qbezier(-1,0)(-0.45,-0.26)(-0.5,-0.866)
\or\qbezier(-1,0)(-0.327,-0.327)(0,-1)
\or\qbezier(-1,0)(-0.179,-0.311)(0.5,-0.866)
\or\qbezier(-1,0)(-0.0536,-0.2)(0.866,-0.5)
\fi
\or\ifcase#2\or\or\or\or\or\or\or
\or\qbezier(-0.866,-0.5)(-0.378,-0.378)(-0.5,-0.866)
\or\qbezier(-0.866,-0.5)(-0.26,-0.45)(0,-1)
\or\qbezier(-0.866,-0.5)(-0.12,-0.446)(0.5,-0.866)
\or\qbezier(-0.866,-0.5)(0,-0.359)(0.866,-0.5)
\fi
\or\ifcase#2\or\or\or\or\or\or\or\or
\or\qbezier(-0.5,-0.866)(-0.138,-0.516)(0,-1)
\or\qbezier(-0.5,-0.866)(0,-0.52)(0.5,-0.866)
\or\qbezier(-0.5,-0.866)(0.12,-0.446)(0.866,-0.5)
\fi
\or\ifcase#2\or\or\or\or\or\or\or\or\or
\or\qbezier(0,-1.)(0.138,-0.516)(0.5,-0.866)
\or\qbezier(0,-1.)(0.26,-0.45)(0.866,-0.5)
\fi
\or\ifcase#2\or\or\or\or\or\or\or\or\or\or
\or\qbezier(0.5,-0.866)(0.378,-0.378)(0.866,-0.5)
\fi\fi\fi\fi}
%
%
\def\FullChord[#1,#2]{
\Endpoint[#1]
\Endpoint[#2]
\Chord[#1,#2]
}
%
%
\def\EndChord[#1,#2]{
\Endpoint[#1]
\Endpoint[#2]
\Chord[#1,#2]
}
%
%
%
%
%
%
\def\Picture#1{
\begin{picture}(2,1)(-1,-0.167)
#1
\end{picture}
}
%
%
\def\DottedChordDiagram[#1,#2]{
\Picture{\DottedCircle \FullChord[#1,#2]}
}
\let\DCD=\DottedChordDiagram
%
%
%
%

\def\SelfInt{\put(-1,-1){\vector(1,1){2}}\put(-1,1){\vector(1,-1){2}}
\put(0,0){\circle*{0.15}}}
\def\Over{\put(-1,-1){\vector(1,1){2}}\put(-1,1){\vector(1,-1){0.75}}\put(0.25,-0.25){\vector(1,-1){0.75}}}
\def\Under{\put(-1,-1){\vector(1,1){0.75}}\put(0.25,0.25){\vector(1,1){0.75}}\put(-1,1){\vector(1,-1){2}}}

\def\Sterm{
\Arc[8]
\Arc[9]
\Arc[10]
\Endpoint[9]
\put(0,-0.16){\circle*{0.15}}
\put(0,-1){\line(0,+1){0.84}}
\put(0,-0.16){\vector(2,1){1}}
\put(0,-0.16){\vector(-2,1){1}}}
\def\Tterm{
\Arc[8]\Arc[9]\Arc[10]
\Endpoint[8]
\Endpoint[10]
}
\def\Tterm{
\Arc[8]
\Arc[9]
\Arc[10]
\Endpoint[8]
\Endpoint[10]
\put(-0.5,-0.866){\vector(-1,2){0.5}}
\put(0.5,-0.866){\vector(1,2){0.5}}
}
\def\Uterm{
\Arc[8]
\Arc[9]
\Arc[10]
\Endpoint[8]
\Endpoint[10]
\put(0.5,-0.866){\vector(-3,2){1.5}}
\put(-0.5,-0.866){\vector(3,2){1.5}}
}

\def\InOut[#1,#2]{
\thinlines
\ifcase#2
\ifcase#1
\qbezier(1,0)(1,0)(0.667,0)
\or\qbezier(1,0)(1,0)(0.206,0.634)
\or\qbezier(1,0)(1,0)(-0.539,0.392)
\or\qbezier(1,0)(1,0)(-0.539,-0.391)
\or\qbezier(1,0)(1,0)(0.206,-0.634)
\or\qbezier(1,0)(1,0)(0.333,0)
\or\qbezier(1,0)(1,0)(-0.167,0.289)
\or\qbezier(1,0)(1,0)(-0.167,-0.289)
\or\qbezier(1,0)(1,0)(0,0)
\fi
\or\ifcase#1
\qbezier(0.866,0.5)(0.866,0.5)(0.667,0)
\or\qbezier(0.866,0.5)(0.866,0.5)(0.206,0.634)
\or\qbezier(0.866,0.5)(0.866,0.5)(-0.539,0.392)
\or\qbezier(0.866,0.5)(0.866,0.5)(-0.539,-0.391)
\or\qbezier(0.866,0.5)(0.866,0.5)(0.206,-0.634)
\or\qbezier(0.866,0.5)(0.866,0.5)(0.333,0)
\or\qbezier(0.866,0.5)(0.866,0.5)(-0.167,0.289)
\or\qbezier(0.866,0.5)(0.866,0.5)(-0.167,-0.289)
\or\qbezier(0.866,0.5)(0.866,0.5)(0,0)
\fi
\or\ifcase#1
\qbezier(0.5,0.866)(0.5,0.866)(0.667,0)
\or\qbezier(0.5,0.866)(0.5,0.866)(0.206,0.634)
\or\qbezier(0.5,0.866)(0.5,0.866)(-0.539,0.392)
\or\qbezier(0.5,0.866)(0.5,0.866)(-0.539,-0.391)
\or\qbezier(0.5,0.866)(0.5,0.866)(0.206,-0.634)
\or\qbezier(0.5,0.866)(0.5,0.866)(0.333,0)
\or\qbezier(0.5,0.866)(0.5,0.866)(-0.167,0.289)
\or\qbezier(0.5,0.866)(0.5,0.866)(-0.167,-0.289)
\or\qbezier(0.5,0.866)(0.5,0.866)(0,0)
\fi
\or\ifcase#1
\qbezier(0,1)(0,1)(0.667,0)
\or\qbezier(0,1)(0,1)(0.206,0.634)
\or\qbezier(0,1)(0,1)(-0.539,0.392)
\or\qbezier(0,1)(0,1)(-0.539,-0.391)
\or\qbezier(0,1)(0,1)(0.206,-0.634)
\or\qbezier(0,1)(0,1)(0.333,0)
\or\qbezier(0,1)(0,1)(-0.167,0.289)
\or\qbezier(0,1)(0,1)(-0.167,-0.289)
\or\qbezier(0,1)(0,1)(0,0)
\fi
\or\ifcase#1
\qbezier(-0.5,0.866)(-0.5,0.866)(0.667,0)
\or\qbezier(-0.5,0.866)(-0.5,0.866)(0.206,0.634)
\or\qbezier(-0.5,0.866)(-0.5,0.866)(-0.539,0.392)
\or\qbezier(-0.5,0.866)(-0.5,0.866)(-0.539,-0.391)
\or\qbezier(-0.5,0.866)(-0.5,0.866)(0.206,-0.634)
\or\qbezier(-0.5,0.866)(-0.5,0.866)(0.333,0)
\or\qbezier(-0.5,0.866)(-0.5,0.866)(-0.167,0.289)
\or\qbezier(-0.5,0.866)(-0.5,0.866)(-0.167,-0.289)
\or\qbezier(-0.5,0.866)(-0.5,0.866)(0,0)
\fi
\or\ifcase#1
\qbezier(-0.866,0.5)(-0.866,0.5)(0.667,0)
\or\qbezier(-0.866,0.5)(-0.866,0.5)(0.206,0.634)
\or\qbezier(-0.866,0.5)(-0.866,0.5)(-0.539,0.392)
\or\qbezier(-0.866,0.5)(-0.866,0.5)(-0.539,-0.391)
\or\qbezier(-0.866,0.5)(-0.866,0.5)(0.206,-0.634)
\or\qbezier(-0.866,0.5)(-0.866,0.5)(0.333,0)
\or\qbezier(-0.866,0.5)(-0.866,0.5)(-0.167,0.289)
\or\qbezier(-0.866,0.5)(-0.866,0.5)(-0.167,-0.289)
\or\qbezier(-0.866,0.5)(-0.866,0.5)(0,0)
\fi
\or\ifcase#1
\qbezier(-1,0)(-1,0)(0.667,0)
\or\qbezier(-1,0)(-1,0)(0.206,0.634)
\or\qbezier(-1,0)(-1,0)(-0.539,0.392)
\or\qbezier(-1,0)(-1,0)(-0.539,-0.391)
\or\qbezier(-1,0)(-1,0)(0.206,-0.634)
\or\qbezier(-1,0)(-1,0)(0.333,0)
\or\qbezier(-1,0)(-1,0)(-0.167,0.289)
\or\qbezier(-1,0)(-1,0)(-0.167,-0.289)
\or\qbezier(-1,0)(-1,0)(0,0)
\fi
\or\ifcase#1
\qbezier(-0.866,-0.5)(-0.866,-0.5)(0.667,0)
\or\qbezier(-0.866,-0.5)(-0.866,-0.5)(0.206,0.634)
\or\qbezier(-0.866,-0.5)(-0.866,-0.5)(-0.539,0.392)
\or\qbezier(-0.866,-0.5)(-0.866,-0.5)(-0.539,-0.391)
\or\qbezier(-0.866,-0.5)(-0.866,-0.5)(0.206,-0.634)
\or\qbezier(-0.866,-0.5)(-0.866,-0.5)(0.333,0)
\or\qbezier(-0.866,-0.5)(-0.866,-0.5)(-0.167,0.289)
\or\qbezier(-0.866,-0.5)(-0.866,-0.5)(-0.167,-0.289)
\or\qbezier(-0.866,-0.5)(-0.866,-0.5)(0,0)
\fi
\or\ifcase#1
\qbezier(-0.5,-0.866)(-0.5,-0.866)(0.667,0)
\or\qbezier(-0.5,-0.866)(-0.5,-0.866)(0.206,0.634)
\or\qbezier(-0.5,-0.866)(-0.5,-0.866)(-0.539,0.392)
\or\qbezier(-0.5,-0.866)(-0.5,-0.866)(-0.539,-0.391)
\or\qbezier(-0.5,-0.866)(-0.5,-0.866)(0.206,-0.634)
\or\qbezier(-0.5,-0.866)(-0.5,-0.866)(0.333,0)
\or\qbezier(-0.5,-0.866)(-0.5,-0.866)(-0.167,0.289)
\or\qbezier(-0.5,-0.866)(-0.5,-0.866)(-0.167,-0.289)
\or\qbezier(-0.5,-0.866)(-0.5,-0.866)(0,0)
\fi
\or\ifcase#1
\qbezier(0,-1)(0,-1)(0.667,0)
\or\qbezier(0,-1)(0,-1)(0.206,0.634)
\or\qbezier(0,-1)(0,-1)(-0.539,0.392)
\or\qbezier(0,-1)(0,-1)(-0.539,-0.391)
\or\qbezier(0,-1)(0,-1)(0.206,-0.634)
\or\qbezier(0,-1)(0,-1)(0.333,0)
\or\qbezier(0,-1)(0,-1)(-0.167,0.289)
\or\qbezier(0,-1)(0,-1)(-0.167,-0.289)
\or\qbezier(0,-1)(0,-1)(0,0)
\fi
\or\ifcase#1
\qbezier(0.5,-0.866)(0.5,-0.866)(0.667,0)
\or\qbezier(0.5,-0.866)(0.5,-0.866)(0.206,0.634)
\or\qbezier(0.5,-0.866)(0.5,-0.866)(-0.539,0.392)
\or\qbezier(0.5,-0.866)(0.5,-0.866)(-0.539,-0.391)
\or\qbezier(0.5,-0.866)(0.5,-0.866)(0.206,-0.634)
\or\qbezier(0.5,-0.866)(0.5,-0.866)(0.333,0)
\or\qbezier(0.5,-0.866)(0.5,-0.866)(-0.167,0.289)
\or\qbezier(0.5,-0.866)(0.5,-0.866)(-0.167,-0.289)
\or\qbezier(0.5,-0.866)(0.5,-0.866)(0,0)
\fi
\or\ifcase#1
\qbezier(0.866,-0.5)(0.866,-0.5)(0.667,0)
\or\qbezier(0.866,-0.5)(0.866,-0.5)(0.206,0.634)
\or\qbezier(0.866,-0.5)(0.866,-0.5)(-0.539,0.392)
\or\qbezier(0.866,-0.5)(0.866,-0.5)(-0.539,-0.391)
\or\qbezier(0.866,-0.5)(0.866,-0.5)(0.206,-0.634)
\or\qbezier(0.866,-0.5)(0.866,-0.5)(0.333,0)
\or\qbezier(0.866,-0.5)(0.866,-0.5)(-0.167,0.289)
\or\qbezier(0.866,-0.5)(0.866,-0.5)(-0.167,-0.289)
\or\qbezier(0.866,-0.5)(0.866,-0.5)(0,0)
\fi
\fi
}

\def\InIn[#1,#2]{
\thinlines
\ifcase#1
\ifcase#2
\qbezier(0.667,0)(0.667,0)(0.667,0)
\or\qbezier(0.667,0)(0.667,0)(0.206,0.634)
\or\qbezier(0.667,0)(0.667,0)(-0.539,0.392)
\or\qbezier(0.667,0)(0.667,0)(-0.539,-0.391)
\or\qbezier(0.667,0)(0.667,0)(0.206,-0.634)
\or\qbezier(0.667,0)(0.667,0)(0.333,0)
\or\qbezier(0.667,0)(0.667,0)(-0.167,0.289)
\or\qbezier(0.667,0)(0.667,0)(-0.167,-0.289)
\or\qbezier(0.667,0)(0.667,0)(0,0)
\fi
\or\ifcase#2
\qbezier(0.206,0.634)(0.206,0.634)(0.667,0)
\or\qbezier(0.206,0.634)(0.206,0.634)(0.206,0.634)
\or\qbezier(0.206,0.634)(0.206,0.634)(-0.539,0.392)
\or\qbezier(0.206,0.634)(0.206,0.634)(-0.539,-0.391)
\or\qbezier(0.206,0.634)(0.206,0.634)(0.206,-0.634)
\or\qbezier(0.206,0.634)(0.206,0.634)(0.333,0)
\or\qbezier(0.206,0.634)(0.206,0.634)(-0.167,0.289)
\or\qbezier(0.206,0.634)(0.206,0.634)(-0.167,-0.289)
\or\qbezier(0.206,0.634)(0.206,0.634)(0,0)
\fi
\or\ifcase#2
\qbezier(-0.539,0.392)(-0.539,0.392)(0.667,0)
\or\qbezier(-0.539,0.392)(-0.539,0.392)(0.206,0.634)
\or\qbezier(-0.539,0.392)(-0.539,0.392)(-0.539,0.392)
\or\qbezier(-0.539,0.392)(-0.539,0.392)(-0.539,-0.391)
\or\qbezier(-0.539,0.392)(-0.539,0.392)(0.206,-0.634)
\or\qbezier(-0.539,0.392)(-0.539,0.392)(0.333,0)
\or\qbezier(-0.539,0.392)(-0.539,0.392)(-0.167,0.289)
\or\qbezier(-0.539,0.392)(-0.539,0.392)(-0.167,-0.289)
\or\qbezier(-0.539,0.392)(-0.539,0.392)(0,0)
\fi
\or\ifcase#2
\qbezier(-0.539,-0.391)(-0.539,-0.391)(0.667,0)
\or\qbezier(-0.539,-0.391)(-0.539,-0.391)(0.206,0.634)
\or\qbezier(-0.539,-0.391)(-0.539,-0.391)(-0.539,0.392)
\or\qbezier(-0.539,-0.391)(-0.539,-0.391)(-0.539,-0.391)
\or\qbezier(-0.539,-0.391)(-0.539,-0.391)(0.206,-0.634)
\or\qbezier(-0.539,-0.391)(-0.539,-0.391)(0.333,0)
\or\qbezier(-0.539,-0.391)(-0.539,-0.391)(-0.167,0.289)
\or\qbezier(-0.539,-0.391)(-0.539,-0.391)(-0.167,-0.289)
\or\qbezier(-0.539,-0.391)(-0.539,-0.391)(0,0)
\fi
\or\ifcase#2
\qbezier(0.206,-0.634)(0.206,-0.634)(0.667,0)
\or\qbezier(0.206,-0.634)(0.206,-0.634)(0.206,0.634)
\or\qbezier(0.206,-0.634)(0.206,-0.634)(-0.539,0.392)
\or\qbezier(0.206,-0.634)(0.206,-0.634)(-0.539,-0.391)
\or\qbezier(0.206,-0.634)(0.206,-0.634)(0.206,-0.634)
\or\qbezier(0.206,-0.634)(0.206,-0.634)(0.333,0)
\or\qbezier(0.206,-0.634)(0.206,-0.634)(-0.167,0.289)
\or\qbezier(0.206,-0.634)(0.206,-0.634)(-0.167,-0.289)
\or\qbezier(0.206,-0.634)(0.206,-0.634)(0,0)
\fi
\or\ifcase#2
\qbezier(0.333,0)(0.333,0)(0.667,0)
\or\qbezier(0.333,0)(0.333,0)(0.206,0.634)
\or\qbezier(0.333,0)(0.333,0)(-0.539,0.392)
\or\qbezier(0.333,0)(0.333,0)(-0.539,-0.391)
\or\qbezier(0.333,0)(0.333,0)(0.206,-0.634)
\or\qbezier(0.333,0)(0.333,0)(0.333,0)
\or\qbezier(0.333,0)(0.333,0)(-0.167,0.289)
\or\qbezier(0.333,0)(0.333,0)(-0.167,-0.289)
\or\qbezier(0.333,0)(0.333,0)(0,0)
\fi
\or\ifcase#2
\qbezier(-0.167,0.289)(-0.167,0.289)(0.667,0)
\or\qbezier(-0.167,0.289)(-0.167,0.289)(0.206,0.634)
\or\qbezier(-0.539,0.392)(-0.167,0.289)(-0.167,0.289)
\or\qbezier(-0.539,-0.391)(-0.167,0.289)(-0.167,0.289)
\or\qbezier(0.206,-0.634)(-0.167,0.289)(-0.167,0.289)
\or\qbezier(0.333,0)(-0.167,0.289)(-0.167,0.289)
\or\qbezier(-0.167,0.289)(-0.167,0.289)(-0.167,0.289)
\or\qbezier(-0.167,-0.289)(-0.167,0.289)(-0.167,0.289)
\or\qbezier(0,0)(-0.167,0.289)(-0.167,0.289)
\fi
\or\ifcase#2
\qbezier(0.667,0)(-0.167,-0.289)(-0.167,-0.289)
\or\qbezier(0.206,0.634)(-0.167,-0.289)(-0.167,-0.289)
\or\qbezier(-0.539,0.392)(-0.167,-0.289)(-0.167,-0.289)
\or\qbezier(-0.539,-0.391)(-0.167,-0.289)(-0.167,-0.289)
\or\qbezier(0.206,-0.634)(-0.167,-0.289)(-0.167,-0.289)
\or\qbezier(0.333,0)(-0.167,-0.289)(-0.167,-0.289)
\or\qbezier(-0.167,0.289)(-0.167,-0.289)(-0.167,-0.289)
\or\qbezier(-0.167,-0.289)(-0.167,-0.289)(-0.167,-0.289)
\or\qbezier(0,0)(-0.167,-0.289)(-0.167,-0.289)
\fi
\or\ifcase#2
\qbezier(0.667,0)(0,0)(0,0)
\or\qbezier(0.206,0.634)(0,0)(0,0)
\or\qbezier(-0.539,0.392)(0,0)(0,0)
\or\qbezier(-0.539,-0.391)(0,0)(0,0)
\or\qbezier(0.206,-0.634)(0,0)(0,0)
\or\qbezier(0.333,0)(0,0)(0,0)
\or\qbezier(-0.167,0.289)(0,0)(0,0)
\or\qbezier(-0.167,-0.289)(0,0)(0,0)
\or\qbezier(0,0)(0,0)(0,0)
\fi
\fi
}

\def\Trythis{
\qbezier(-1,0)(-1,0)(0.866,-0.5)
\put(-0.0536,-0.2){\circle*{0.15}}
}


\def\DoDash[#1,#2]{
\ifcase#1
\ifcase#2
\qbezier[20](0.5,0)(0,-1.25)(-0.5,0)
\or\qbezier[20](0.5,0)(-0.25,-1.25)(-1,-0.5)
\or\qbezier[20](0.5,0)(+0.25,-1.25)(-0.25,-0.5)
\or\qbezier[20](0.5,0)(0,-2.0)(-1.5,-1)
\or\qbezier[20](0.5,0)(0,-2.0)(-0.75,-1)
\or\qbezier[20](0.5,0)(0.25,-2.0)(-0.5,-1)
\or\qbezier[20](0.5,0)(0.25,-2.0)(0,-1)
\fi
\or\ifcase#2
\qbezier[20](0.25,-0.5)(-0.25,-1.25)(-0.5,0)
\or\qbezier[20](0.25,-0.5)(-0.25,-1.25)(-1,-0.5)
\or\qbezier[20](0.25,-0.5)(0,-1.25)(-0.25,-0.5)
\or\qbezier[20](0.25,-0.5)(0,-2.0)(-1.5,-1)
\or\qbezier[20](0.25,-0.5)(0,-2.0)(-0.75,-1)
\or\qbezier[20](0.25,-0.5)(0,-2.0)(-0.5,-1)
\or\qbezier[20](0.25,-0.5)(0.25,-2.0)(0,-1)
\fi\or\ifcase#2
\qbezier[20](1,-0.5)(0,-1.5)(-0.5,0)
\or\qbezier[20](1,-0.5)(0,-1.5)(-1,-0.5)
\or\qbezier[20](1,-0.5)(0,-1.5)(-0.25,-0.5)
\or\qbezier[20](1,-0.5)(0,-2.0)(-1.5,-1)
\or\qbezier[20](1,-0.5)(0,-2.0)(-0.75,-1)
\or\qbezier[20](1,-0.5)(0,-2.0)(-0.5,-1)
\or\qbezier[20](1,-0.5)(0,-2.0)(0,-1)
\fi
\fi
}

\def\DoA{
\put(0,1){\line(0,-1){0.5}}
}
\def\DoB{
\put(-0.5,0){\line(1,1){0.5}}
}\
\def\DoC{
\put(0,0.5){\line(1,-1){0.5}}
}
\def\DoD{
\put(-1,-0.5){\line(1,1){0.5}}
}
\def\DoE{
\put(-0.5,0){\line(1,-2){0.25}}
}
\def\DoF{
\put(0.25,-0.5){\line(1,2){0.25}}
}
\def\DoG{
\put(0.5,0){\line(1,-1){0.5}}
}
\def\DoH{
\put(-1.5,-1){\line(1,1){0.5}}
}
\def\DoI{
\put(-1,-0.5){\line(1,-2){0.25}}
}
\def\DoJ{
\put(-0.5,-1){\line(1,2){0.25}}
}
\def\DoK{
\put(-0.25,-0.5){\line(1,-2){0.25}}
}


\begin{titlepage}
\hskip 1cm

\begin{flushright}
QMW-PH-95-32\\
UM-P-95-112 \\
q-alg/9511024\\
\vspace{12pt}
First Edition: 5 November 1995\\
This Edition: 19 September 1996
\end{flushright}

\hskip 1.5cm
\vfill
\begin{center}
{\LARGE Cabling the Vassiliev Invariants
}\vfill
{\large A. Kricker,$^1$
B. Spence$^2$ and I.Aitchison$^3$\\
\vspace{6pt}
}
\end{center}

\vfill
\begin{center}
{\bf Abstract}
\end{center}
{\small
We characterise the cabling operations on the weight systems of finite type
knot invariants. 
The eigenvectors and eigenvalues of this family of
operations are described. The canonical deframing projection for these knot
invariants is described over the cable eigenbasis. The action of immanent
weight systems on general Feynman diagrams is considered, and
the highest eigenvalue cabling eigenvectors are shown to be dual to
the immanent weight systems. 
Using these results, we prove a recent conjecture of Bar-Natan and Garoufalidis
on cablings of weight systems.
}

 \vfill \hrule width 3.cm
{\footnotesize
\noindent
 $^1$ School of Mathematical Sciences, University of Melbourne, Parkville
3052 Australia.
\\
\hphantom{$^1$}
Email: krick@mozart.ph.unimelb.edu.au\\
\noindent$^2$ 
Department of Physics, Queen Mary and Westfield College, London E1 4NS,
UK.
\\
\hphantom{$^1$}
 Email:
B.Spence@qmw.ac.uk.
\\
\noindent$^3$ 
School of Mathematical Sciences, University of Melbourne, Parkville 3052
Australia.
\\
\hphantom{$^1$}
 Email:
iain@mundoe.maths.mu.oz.au
}

\end{titlepage}


\section{Introduction.}

The theory of knot invariants of finite order (or Vassiliev invariants \cite{Vas})
is today well established. In particular, the essential combinatorial
structure of the constructions has been elucidated by Birman and Lin
\cite{BL}, and 
Bar-Natan \cite{BN1}. Kontsevich showed that, up to knot invariants
of lesser order, these combinatorial representations are faithful \cite{Kon}.
Several important questions remain.
First amongst these is  the question of whether 
all knot invariants of finite order are realisable as finite linear
combinations of the coefficients of Lie algebraic (quantum group) knot
invariants. Another question is whether a better  understanding of
 Vassiliev knot
invariants would follow
from a more traditional algebraic-topological methodology.

In a 1994 paper \cite{BNG}, Bar-Natan and Garoufalidis presented the following 
theorem, originally conjectured by Melvin and Morton \cite{MM}. Note that their proof followed a path-integral demonstration by Rozansky \cite{Roz}. 

 We denote by $\EJ(K)(\hb)$ 
the $U_q(sl(2))$ invariant evaluated in a representation of dimension $\lambda
+ 1$ at $q\ =\ e^{\hb}$ and we will denote the Alexander polynomial as
$A(K)(z)$.

\begin{theorem}[\cite{BNG}]
Expanding $\EJ/(\lambda + 1)$ in powers of $\lambda$ and $\hb$,
\bq
\frac{\EJ(K)(\hb)}{\lambda + 1}\ =\ \sum_{j,m\geq 0} b_{jm}(K)\lambda^j \hb^m,
\eq
we have,
\begin{enumerate}
\item{$b_{jm}(K) = 0$ if $j>m$.}
\item{Define \[ JJ(K)(\hb) = \sum_{m=0}^\infty b_{mm}(K) \hb^m. \]

	Then,
\bq
JJ(K)(\hb)\,\frac{\hb}{e^{\frac{\hb}{2}}-e^{-\frac{\hb}{2}}}\,A(K)(e^{\hb}) = 1.
\eq
}
\end{enumerate}
\rtb
\end{theorem}

Thus Bar-Natan and Garoufalidis showed that one can reconstruct the
Alexander-Conway polynomial from the highest order terms in the
extended Jones polynomial, when that polynomial is expanded in the
dimension of the representation.
They also extended this result to any semi-simple
Lie algebra.

 There is a certain philosophy motivated by extensive numerical calculations
that the Lie algebraic invariants should span
the set of finite order invariants. Thus their result led
 Bar-Natan and Garoufalidis to
conjecture that there was in some natural sense a highest-order term in an
arbitrary Vassiliev invariant expressible over the algebra of coefficients of
the Alexander-Conway polynomial. They 
conjectured
 that the order of the natural cabling operation would be the intrinsic
variable in which one could expand weight systems, replacing the
dimension of the representation.

We will prove the following formulation of this conjecture.

\begin{theorem}

Let $W_m$ be a weight system of order m (a weight system on m-chorded
diagrams, $A_m$). Denote for the deframed n-th cable of this weight system,
\bq
\widehat{\psi^{n*}W_m}\ =\ W_m\circ\psi^n_m\circ \phi_m.
\eq
In the above $\phi_m$ is a deframing projector (which we shall discuss below).
\begin{enumerate}
\item{ As a function of n, $\widehat{\psi^{n*}W_m}$ is a polynomial in n, of
highest order $n^m$.}
\item{ The coefficient of $n^{m}$ in the polynomial is equal to a linear
combination of immanent weight systems. }
\end{enumerate}
\rtb
\end{theorem}

Immanent weight systems were written down in \cite{BNG}. They are a means of
calculating invariants of chord diagrams from their intersection matrices.
These authors also showed that the algebra of immanent weight systems is the algebra 
of the weight systems for the coefficients of
the Alexander polynomial. That is to say, the subspace of weight systems coming
from sums of products of coefficients of the Alexander-Conway polynomial is
the same as the subspace coming from immanent weight systems.
Details of this can be found in the discussion in Section 6 of \cite{BNG}.

In the following, we will begin in section 2 by reviewing briefly
some essential concepts - the notions of the chord diagram
algebra, weight systems and Feynman diagrams.
In section 3, we introduce cabling, and construct eigenvectors
of the cabling operation (our results in this section partially answer
Bar-Natan's Problem 7.4 \cite{BN1}, providing a topological understanding
of the alternative grading on the chord diagram algebra).

These eigenvectors have in fact already been written down
in another context, relating to the enumeration of
 primitive vectors for this algebra.
In section 4, we consider the deframing operation on the chord diagram
algebra. The action of deframing finds a simple expression on our collection
of eigenvectors. We use this to prove  the first part of Theorem $2$, 
enumerating
the eigenvectors that are in the deframing invariant subspace.
In section 5, we recall immanent weight systems, and their relation with the
Alexander-Conway weight system. In this section we show that 
they form an orthogonal dual  basis on
precisely the eigenvectors 
 within the highest eigenvalue subspace,
and are zero on all others. This concludes the proof of Theorem 2.
Section 6 contains the analysis of the cycle decomposition
sums of FDs: this facilitates the proof of many statements made
in section 5. We finish
by discussing some further directions of research suggested by these results.


\section{ The Chord Diagram Algebra. }

In what follows we employ a representative field of characteristic 0, the
complex numbers 
$\ccc$ (unless stated otherwise). 
We shall consider framed knots: when we say knot invariant, it may 
be framing dependent.
We begin with knot invariants of
finite order. To facilitate this definition we first extend any
knot invariant $V$ to an invariant of knots with self-intersections
(with blackboard framing): 

\

\begin{equation}
V(\ \begin{picture}(2,1)(-1,-0.167) \SelfInt \end{picture}\ )\ =\ 
V(\ \begin{picture}(2,1)(-1,-0.167) \Under \end{picture}\ )\ -\
V(\ \begin{picture}(2,1)(-1,-0.167) \Over \end{picture}\ )\ .\
\end{equation}

\

\begin{defn} 

A knot invariant is {\bf of finite order n}, if it 
vanishes on knots with more than $n$ self-intersections.
\rtb
\end{defn}

Denote by $\hat{\aaa}_m$ the finite-dimensional vector space spanned by
 $\ccc$-linear sums of chord diagrams with $m$ chords. The loop is 
oriented; in diagrams we shall assume counterclockwise. For
example, 

\begin{equation}
\Picture{
\FullCircle
\FullChord[0,8]
\FullChord[2,6]
\FullChord[4,10]}
\ -\ 2\ 
\Picture{
\FullCircle
\FullChord[0,6]
\FullChord[2,8]
\FullChord[4,10]
}
\ +\pi \
\Picture{
\FullCircle
\FullChord[0,2]
\FullChord[4,6]
\FullChord[8,10]
}\ \ci\ \ddd_3.
\end{equation}

\

Chord diagrams, in a natural way, represent the order in which
intersections are encountered in a singular knot. An invariant
of finite type $n$, provides a well-defined dual vector to $\hat{\aaa}_n$. 
Namely, such an invariant will only observe the order in which 
self-intersections are encountered, and not extra ``knotting" information. 
The invariant is then linearly extended to linear combinations of diagrams. 

There are certain relationships that the resulting functional satisfies if it
is obtained in this way: the 4-T relations, which
can be understood most easily from a three-dimensional picture \cite{BN1}. 
Generalising:
\begin{defn}

A $\ccc$-valued {\bf weight system} of degree m, is a linear functional
$W: \ddd_m \longrightarrow\ \ccc$ satisfying the 4-T relations.

\

{\bf \underline{4-T}}.

\begin{equation}
W\left(
\Picture{\DottedCircle\FullChord[1,8]\Arc[2]\FullChord[5,9]
\Arc[1]\Arc[8]\Arc[5]\Arc[9]}
\right) -
W\left(
\Picture{\DottedCircle\FullChord[1,9]\Arc[2]\FullChord[5,8]
\Arc[1]\Arc[9]\Arc[5]\Arc[8]}
\right) +
W\left(
\Picture{\DottedCircle\FullChord[2,5]\Arc[8]\FullChord[1,9]
\Arc[2]\Arc[5]\Arc[1]\Arc[9]}
\right) -
W\left(
\Picture{\DottedCircle\FullChord[1,5]\Arc[8]\FullChord[2,9]
\Arc[1]\Arc[5]\Arc[2]\Arc[9]}
\right) = 0.
\end{equation}

\rtb
\end{defn}

With regard to pictures of chord diagrams in this paper, sections of the 
outer circle of a chord diagram which are dotted denote parts of the
outer circle where further chords or lines not shown in the diagram may end.
Sections of the outer circle which are full lines show regions where
all allowed terminal points of chords or lines are shown.  In what follows
we shall refer to the outer loop as the {\bf Wilson loop} of the diagram.

We denote the degree $m$ weight system that one constructs from a given knot
invariant of 
finite order $m$, $V$, by $W_m[V]$. The vector space dual to the space of 
weight systems at 
degree $m$ 
is denoted by $\aaa_m$ (i.e. $\ddd_m$ quotiented by 4-T expressions).

It is convenient to express frequently recurring
 linear combinations of diagrams
in $\aaa_m$ using another notation. The more general
diagrams are 
referred to as {\bf Feynman Diagrams }(FDs), and allow internal
 trivalent
vertices, 
where the incoming edges to the vertex are assigned a cyclic ordering. 

The following definitions will be useful.
\begin{defn}
The {\bf graph} of a FD is the graph that remains when the Wilson loop
is removed. It has univalent and trivalent vertices, and the incoming
edges to a trivalent vertex are cyclically ordered. A {\bf leg} is a
neighbourhood in the graph of some univalent vertex.
Denote by $\ddd^t_m$ the vector space whose basis
is formed by FDs whose graphs contain $2m$ vertices.
\rtb
\end{defn}

Let the STU vectors of $\ddd^t_m$ be the following (where only a part of the
graph is shown): 

\begin{equation}
\Picture{
\DottedCircle
\InOut[3,7]
\InIn[3,2]
\InIn[3,4]
\Arc[7]\Endpoint[7]
}\ \ -\ \
\Picture{
\DottedCircle
\InOut[2,7]
\InOut[4,8]
\Arc[8]\Arc[7]
\Endpoint[8]\Endpoint[7]
}\ \ +\ \
\Picture{
\DottedCircle
\InOut[2,8]
\InOut[4,7]
\Arc[8]\Arc[7]
\Endpoint[8]\Endpoint[7]
}
\end{equation}

\

\

\begin{defn}
Denote by $\aaa^t_m$ the quotient of $\ddd^t_m$ by the subspace spanned
by the STU
vectors. According to \cite{BN1}, $\aaa_m\equiv \aaa^t_m$. Hereafter
we shall only refer to $\aaa_m$.
\rtb
\end{defn}

A FD with
$2m$ vertices (internal and external) {\it resolves} to an $m$-chorded diagram.
For example, the following diagram in $\aaa_3$ can be expanded as

\begin{eqnarray}
\left.
\Picture{
\FullCircle
\InIn[7,6]
\InOut[7,8]
\InOut[7,10]
\InOut[6,2]
\InOut[6,4]
\Endpoint[8]\Endpoint[10]\Endpoint[2]\Endpoint[4]
}\right. & = &
\left.
\Picture{
\FullCircle
\InOut[6,2]
\InOut[6,4]
\InOut[6,8]
\FullChord[9,10]
\Endpoint[2]\Endpoint[4]\Endpoint[8]
}\ \ -\
\Picture{
\FullCircle
\InOut[6,2]
\InOut[6,4]
\InOut[6,9]
\FullChord[8,10]
\Endpoint[2]
\Endpoint[4]
\Endpoint[9]
}\right. \nonumber \\
& & \nonumber \\
& & \nonumber \\
 & = &
\left.
\Picture{
\FullCircle
\FullChord[4,7]
\FullChord[2,8]
\FullChord[9,11]
}\ -\ 
\Picture{
\FullCircle
\FullChord[2,7]
\FullChord[4,8]
\FullChord[9,11]
}\ -\
\Picture{
\FullCircle
\FullChord[4,8]
\FullChord[2,9]
\FullChord[7,11]
}\ +\
\Picture{
\FullCircle
\FullChord[7,11]
\FullChord[4,9]
\FullChord[2,8]
}
\right.\ .
\end{eqnarray}

\

We equip $\aaa_m$ with a product, being the connect-sum of Wilson loops
at some choice of point, 4-T relations being required in order that this product be well-defined. For example:  

\begin{equation}
\Picture{
\FullCircle
\InOut[8,8]
\InOut[8,10]
\InOut[8,3]
\Endpoint[8]\Endpoint[10]\Endpoint[3]
}\ .\
\Picture{
\FullCircle
\FullChord[8,2]
\FullChord[10,4]
}
\ \ =\
\ \Picture{
\FullCircle
\InOut[3,4]
\InOut[3,7]
\InOut[3,9]
\FullChord[1,10]
\FullChord[11,2]
\Endpoint[4]
\Endpoint[7]
\Endpoint[9]
}\ \ \ci\ \ \aaa_4.
\end{equation}

\

This forms a graded algebra with identity element the empty Wilson loop.
Furthermore, there is a co-product homomorphism $\Delta:
A\longrightarrow A\otimes A$. It is defined 
as a sum over all ways of partitioning the chords of a diagram into two sets,
and 
separating them accordingly.
For example:

\begin{eqnarray}
\Delta\left( 
\Picture{
\FullCircle
\FullChord[3,7]
\FullChord[2,8]
\FullChord[5,11]
}\right) & = & 
\ \ \, \Picture{
\FullCircle
\FullChord[3,7]
\FullChord[2,8]
\FullChord[5,11]
}\otimes
\ \, \Picture{
\FullCircle}
\ \ +\ \ 2\, 
\Picture{
\FullCircle
\FullChord[2,8]
\FullChord[5,11]
}\otimes
\ \, \Picture{
\FullCircle
\FullChord[3,9]
}\ \ \ \
 \nonumber\\
& & \nonumber \\
& & \nonumber \\
& & 
+ \ \, \Picture{
\FullCircle
\FullChord[3,7]
\FullChord[2,8]
}\otimes 
\ \, \Picture{
\FullCircle
\FullChord[3,9]
}\
+\ \ \, \Picture{
\FullCircle
\FullChord[3,9]
}
\otimes
\ \, \Picture{
\FullCircle
\FullChord[3,7]
\FullChord[2,8]
}\ \ \ \
 \nonumber\\
& & \nonumber \\
& & \nonumber \\
& & 
+ 2\, \Picture{
\FullCircle
\FullChord[3,9]
}\otimes
\ \, \Picture{
\FullCircle
\FullChord[5,11]
\FullChord[2,8]
}\ \ +\ \ \, \Picture{
\FullCircle}
\otimes
\ \, \Picture{
\FullCircle
\FullChord[2,8]
\FullChord[3,7]
\FullChord[5,11]
}.
\nonumber \\
& & \nonumber \\
& &
\end{eqnarray}

These operations are well-defined and satisfy
the axioms of a (commutative and co-commutative)
Hopf algebra. By the structure theory of Hopf algebras,
such an algebra is generated by its primitive elements,
namely those $v\,\ci\, \aaa_{*}$ such
that
\bq
\Delta(v)\ =\ v\,\otimes\,1\ +\ 1\,\otimes\,v.
\eq
In the next section we shall describe the primitive subspace in terms of 
eigenvectors of certain cabling operations.


\section{Cabling}

Cabling is a natural 
topological operation on framed knots. Moreover, when composed with knot
invariants of finite order, it leads to an interesting collection of linear
transformations on the weight systems that characterise those knot
invariants. 

\begin{defn}
Consider a knot $z^\mu:(0,2\pi]\longrightarrow \rrr^3$. Observe
that the normal plane sufficiently close to a point on the knot can be
parameterised as a domain in the complex plane, with the knot intersecting at
the origin. With this, take a  framing of the knot $\epsilon: (0,2\pi]
\longrightarrow \ccc$, and the embedding of the normal plane in $\rrr^3,\
e^\mu: \ccc \times (0,2\pi] \longrightarrow \rrr^3$.

\

The {\bf nth-connected  cabling} of $z^{\mu}$ is the knot given by
\bq
(\psi^n z^\mu)(t)\ =\ z^{\mu}(nt) + e^\mu(nt,e^{2\pi i t}\epsilon(nt)),
\eq
(with the understanding that angles are identified mod $2\pi$).
\rtb
\end{defn}

Operating inside $V$, a (framing independent) knot invariant of 
finite order $m$, this yields 
$V\circ \psi^n$, a framing
dependent invariant also of finite order $m$.
We wish to describe 
an operator $\psi^n_m$
satisfying 
\bq
W_m[V\circ \psi^n]\ =\ W_m[V]\circ \psi^n_m\ .
\label{cabledefn}
\eq

\begin{defn}

Denote by $\psi^n_m$ the linear transformation, $\psi^n_m: A_m
\longrightarrow A_m$ defined as follows. Take the nth cyclic cover of $S^1$,
which is also $S^1$. Sum over all ways of ``lifting'' the ends of chords to
the different covers.
For example, 
\begin{eqnarray}
\psi^2_2 \left( 
\Picture{
\FullCircle
\FullChord[2,8]
\FullChord[5,11]
}\right) & = &
\Picture{
\CabCirc
\CabE
\CabF
}+
\Picture{
\CabCirc
\CabE
\CabA
}+
\Picture{
\CabCirc
\CabE
\CabB
}+
\Picture{
\CabCirc
\CabF
\CabC
}+\ \ldots
\nonumber \\
&  & \nonumber \\
& = & 
8\,\Picture{
\FullCircle
\FullChord[2,8]
\FullChord[4,10]
}\ +\ 
8\, \Picture{
\FullCircle
\FullChord[2,10]
\FullChord[4,8]
} \nonumber \\
& & \end{eqnarray}

\
This operation satisfies equation  (\ref{cabledefn}) \cite{BN1}. 
\rtb
\end{defn}

The cabling transformation is well defined on $\aaa_{*}$, with the
4-T relations being mapped into themselves.
To see this,  we describe 
cabling on FDs with internal trivalent vertices.

Consider first 
the effect of cabling on a single trivalent vertex:   

\begin{eqnarray}
\psi^2\left(
\Picture{
\put(-1,-1){\vector(1,0){2}}
\put(0,-1){\circle*{0.15}}
\put(0,0){\line(0,-1){1}}
\put(-1,0.5){\line(2,-1){1}}
\put(0,0){\line(2,1){1}}
}\right)
 & = & \psi^2
\left(
\Picture{
\put(-1,-1){\vector(1,0){2}}
\put(-1,0.5){\line(1,-2){0.75}}
\put(-0.25,-1){\circle*{0.15}}
\put(0.25,-1){\circle*{0.15}}
\put(0.25,-1){\line(1,2){0.75}}
}\ -\
\Picture{
\put(-1,-1){\vector(1,0){2}}
\put(-0.25,-1){\circle*{0.15}}
\put(0.25,-1){\circle*{0.15}}
\put(-1,0.5){\line(5,-6){1.25}}
\put(-0.25,-1){\line(5,6){1.25}}
}\right) \nonumber \\
& & \nonumber \\
& = &
\Picture{
\put(-1,-1){\vector(1,0){2}}
\put(-1,-0.5){\vector(1,0){2}}
\put(-0.25,-1){\circle*{0.15}}
\put(0.25,-1){\circle*{0.15}}
\put(-1,0.5){\line(1,-2){0.75}}
\put(0.25,-1){\line(1,2){0.75}}
}\ +\
\Picture{
\put(-1,-1){\vector(1,0){2}}
\put(-1,-0.5){\vector(1,0){2}}
\put(-0.25,-1){\circle*{0.15}}
\put(0.25,-0.5){\circle*{0.15}}
\put(-1,0.5){\line(1,-2){0.75}}
\put(0.25,-0.5){\line(3,4){0.75}}
}\ +\
\Picture{
\put(-1,-1){\vector(1,0){2}}
\put(-1,-0.5){\vector(1,0){2}}
\put(-0.25,-0.5){\circle*{0.15}}
\put(0.25,-1){\circle*{0.15}}
\put(-1,0.5){\line(3,-4){0.75}}
\put(0.25,-1){\line(1,2){0.75}}
}\ +\
\Picture{
\put(-1,-1){\vector(1,0){2}}
\put(-1,-0.5){\vector(1,0){2}}
\put(-0.25,-0.5){\circle*{0.15}}
\put(0.25,-0.5){\circle*{0.15}}
\put(-1,0.5){\line(3,-4){0.75}}
\put(0.25,-0.5){\line(3,4){0.75}}
} \nonumber \\
& & \nonumber \\
& & 
-\Picture{
\put(-1,-1){\vector(1,0){2}}
\put(-1,-0.5){\vector(1,0){2}}
\put(-0.25,-1){\circle*{0.15}}
\put(0.25,-1){\circle*{0.15}}
\put(-1,0.5){\line(5,-6){1.25}}
\put(-0.25,-1){\line(5,6){1.25}}
}\ -\
\Picture{
\put(-1,-1){\vector(1,0){2}}
\put(-1,-0.5){\vector(1,0){2}}
\put(-0.25,-1){\circle*{0.15}}
\put(0.25,-0.5){\circle*{0.15}}
\put(-1,0.5){\line(1,-2){0.75}}
\put(0.25,-0.5){\line(3,4){0.75}}
}\ -\
\Picture{
\put(-1,-1){\vector(1,0){2}}
\put(-1,-0.5){\vector(1,0){2}}
\put(-0.25,-0.5){\circle*{0.15}}
\put(0.25,-1){\circle*{0.15}}
\put(-1,0.5){\line(3,-4){0.75}}
\put(0.25,-1){\line(1,2){0.75}}
}\ -\
\Picture{
\put(-1,-1){\vector(1,0){2}}
\put(-1,-0.5){\vector(1,0){2}}
\put(-0.25,-0.5){\circle*{0.15}}
\put(0.25,-0.5){\circle*{0.15}}
\put(-1,0.5){\line(5,-4){1.25}}
\put(-0.25,-0.5){\line(5,4){1.25}}
} \nonumber \\
& & \nonumber \\
&=&
\Picture{
\put(-1,-1){\vector(1,0){2}}
\put(-1,-0.5){\vector(1,0){2}}
\put(-0.25,-1){\circle*{0.15}}
\put(0.25,-1){\circle*{0.15}}
\put(-1,0.5){\line(1,-2){0.75}}
\put(0.25,-1){\line(1,2){0.75}}
}\ -\
\Picture{
\put(-1,-1){\vector(1,0){2}}
\put(-1,-0.5){\vector(1,0){2}}
\put(-0.25,-1){\circle*{0.15}}
\put(0.25,-1){\circle*{0.15}}
\put(-1,0.5){\line(5,-6){1.25}}
\put(-0.25,-1){\line(5,6){1.25}}
}\ +\
\Picture{
\put(-1,-1){\vector(1,0){2}}
\put(-1,-0.5){\vector(1,0){2}}
\put(-0.25,-0.5){\circle*{0.15}}
\put(0.25,-0.5){\circle*{0.15}}
\put(-1,0.5){\line(3,-4){0.75}}
\put(0.25,-0.5){\line(3,4){0.75}}
}\ -\
\Picture{
\put(-1,-1){\vector(1,0){2}}
\put(-1,-0.5){\vector(1,0){2}}
\put(-0.25,-0.5){\circle*{0.15}}
\put(0.25,-0.5){\circle*{0.15}}
\put(-1,0.5){\line(5,-4){1.25}}
\put(-0.25,-0.5){\line(5,4){1.25}}
} \nonumber \\
& & \nonumber \\
& & 
\end{eqnarray}

It is clear from this picture that in order to cable a diagram whose graph has
a single trivalent vertex one  
\lq sums over 
lifts' of the legs of $v$.
By induction this understanding extends to diagrams with an arbitrary 
number of internal vertices.

Denote by $S_p$ the group of permutations of $p$ letters.
There exists a natural action of $S_p$ on \lq line' chord diagrams, with $p$ 
legs (this action depends 
on the point where one breaks the Wilson loop).
For example, we may break the following three chord diagram.
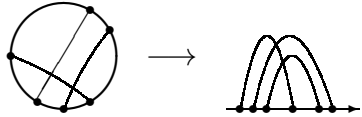
\[
\Picture{
\FullCircle
\FullChord[6,10]
\FullChord[9,1]
\FullChord[8,2]
}\ \ \lra\ \ 
\Picture{
\put(-1,-1){\vector(1,0){2.5}}
\qbezier(-0.75,-1)(-0.25,1.75)(0.25,-1)
\qbezier(-0.5,-1)(0.1,1.75)(1.0,-1)
\qbezier(-0.25,-1)(0.2,1.0)(0.75,-1)
\put(-0.75,-1){\circle*{0.15}}
\put(0.25,-1){\circle*{0.15}}
\put(-0.5,-1){\circle*{0.15}}
\put(1,-1){\circle*{0.15}}
\put(-0.25,-1){\circle*{0.15}}
\put(0.75,-1){\circle*{0.15}}
}
\]

\

Elements of $S_p$ then act by permuting the order of the location of the
univalent vertices (the legs) on the Wilson loop. 
To construct a cabling eigenvector, let us start with some choice of a Feynman
diagram, $v\in \aaa_n$, with $p$ legs.  Construct the vector
$Sym_{v}\in \aaa_n$ as follows. 

\begin{defn}

\bq
Sym_{v}\ =\ \sum_{\sigma \ci S_p} \sigma(v).
\eq
\rtb
\end{defn}

(Note that this is independent of the choice of break because we have summed over all
permutations.)
This vector is the same for initial FDs related by a permutation of their legs.
The graphs of two such FDs have the same set of connected
components when the Wilson loop is removed. We record these vectors in this
way: the collection
of graphs with \lq internal', trivalent vertices and \lq external', univalent
vertices. In fact, they have been written down previously, in another context,
by Kontsevich \cite{Kon}. There they were used to describe the primitive
vectors of the chord diagram algebra, and were dubbed \lq Chinese
Characters' by Bar-Natan, \cite{BN1}. Here we shall refer to them as
symmetrised Feynman diagrams (SFDs). Examples are
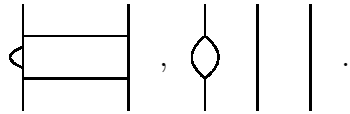
\[
\Picture{
\put(-1,1){\line(0,-1){2}}
\put(1,1){\line(0,-1){2}}
\put(-1,0.4){\line(1,0){2}}
\put(-1,-0.4){\line(1,0){2}}
\qbezier(-1,0.2)(-1.5,0)(-1,-0.2)
}\ \ \ ,\ \ \
\Picture{
\put(-1,1){\line(0,-1){0.6}}
\put(-1,-0.4){\line(0,-1){0.6}}
\qbezier(-1,0.4)(-1.5,0)(-1,-0.4)
\qbezier(-1,0.4)(-0.5,0)(-1,-0.4)
\put(0,1){\line(0,-1){2}}
\put(1,1){\line(0,-1){2}}
}\ \ \ .\ \ 
\]

\

Keeping in mind their representation as sums over ways of ordering the 
univalent vertices
of the graph on a Wilson loop, and the STU relations, 
there are further identities of importance in this space. In the
diagrams below, we have extracted one part of a SFD, 
further connections from the univalent vertices to other parts of the diagram 
are possible. 

\

{\bf Antisymmetry.}
\begin{equation}
\Picture{
\put(0,0){\line(0,1){1}}
\put(0,0){\line(1,-1){1}}
\put(0,0){\line(-1,-1){1}}
}\ \ =\ \
\ \ -\ \Picture{
\put(0,0){\line(0,1){1}}
\qbezier(0,0)(-1,-0.5)(1,-1)
\qbezier(0,0)(1,-0.5)(-1,-1)
\put(1,-1){\line(2,-1){0.1}}
\put(-1,-1){\line(-2,-1){0.1}}
}
\end{equation}

\

{\bf IHX relation.}
\begin{equation}
\Picture{
\put(-1,1){\line(1,0){2}}
\put(-1,-1){\line(1,0){2}}
\put(-1,1){\line(-1,0){0.1}}
\put(-1,-1){\line(-1,0){0.1}}
\put(0,1){\line(0,-1){2}}
}
\ \ =\ \ 
\Picture{
\qbezier(-1,1)(-0.5,1)(-0.5,0)
\qbezier(-1,-1)(-0.5,-1)(-0.5,0)
\qbezier(1,1)(0.5,1)(0.5,0)
\qbezier(1,-1)(0.5,-1)(0.5,0)
\put(-0.5,0){\line(1,0){1}}
\put(-1,1){\line(-1,0){0.1}}
\put(-1,-1){\line(-1,0){0.1}}
\put(1,1){\line(1,0){0.1}}
\put(1,-1){\line(1,0){0.1}}
}\ \ -\ \
\Picture{
\qbezier(-1,1)(-0.5,1)(-0.25,0)
\qbezier(-0.25,0)(0,-1)(1,-1)
\qbezier(1,1)(0.5,1)(0.25,0)
\qbezier(-1,-1)(0,-1)(0.25,0)
\put(-0.25,0){\line(1,0){0.5}}
\put(-1,1){\line(-1,0){0.01}}
\put(-1,-1){\line(-1,0){0.01}}
\put(1,1){\line(1,0){0.01}}
\put(1,-1){\line(1,0){0.01}}
}
\end{equation}
\rtb

\

\begin{theorem}

If the graph of some diagram $v$ has $p$ univalent vertices, then \[
 \psi^{n}_{m}\, Sym_{v}\ =\ n^p\, Sym_{v}. \]

\end{theorem}
{\bf \underline{Proof}.}

The action of cabling on any Feynman diagram is expressible as a sum over
actions by elements of $S_p$, once a labelling of the legs has been 
chosen. Recall that $Sym_{v}$ is a sum
over all possible orderings of the external edges of a FD. Thus the action of
any $\sigma\, \ci\, S_p$ just returns $Sym_{v}$. To decide the eigenvalue
it remains to count the number of permutations by which
 a cabling is expressed:
every leg can  be raised to one of $n$ possible covers, and there are
p legs, so $\psi^n_m$ is expressed on a FD with $p$ legs as a sum of $n^p$
permutation actions. Thus $n^p$ is the eigenvalue.
\rtb

Hence we have a useful collection of eigenvectors of the cabling operation. 
In fact, we have a {\it diagonalisation}. In {\cite{BN1}} it is shown how
the difference between any FD and a permutation of its legs is
expressible as a sum of FDs with fewer legs. It follows
that any FD can be expressed as a sum of these totally symmetrised vectors.


\section{Deframing}

We have seen that cabling is an operation on framed knots. We require a
framing to choose a particular cabling, and so cabling a knot invariant, even
if it previously was framing independent, introduces a framing dependence.

Framing independence of a Vassiliev knot invariant translates into
an additional set of relations on the associated weight system, the 1-T
relations.
Diagrammatically,

\begin{equation}
W[V] \left(
\Picture{
\DottedCircle
\FullChord[6,7]
\Arc[6]
\Arc[7]}
\right)
\ =\  
V\left(
\Picture{
\qbezier(-0.8,0.9)(-0.8,0.45)(0,0)
\qbezier(-0.8,-0.9)(-0.8,-0.45)(0,0)
\put(0,0){\circle*{0.15}}
\qbezier(0,0)(1,0.5)(1,0)
\qbezier(0,0)(1,-0.5)(1,0)}
\right)
\  =\ 
V\left(
\Picture{
\qbezier(-0.8,0.9)(-0.8,0.45)(-0.2,0.1)
\qbezier(0.2,-0.10)(1,-0.5)(1,0)
\qbezier(-0.8,-0.9)(-0.8,-0.45)(0,0)
\qbezier(0,0)(1,0.5)(1,0)
}\right)\ -\
V\left(
\Picture{
\qbezier(-0.8,0.9)(-0.8,0.45)(0,0)
\qbezier(0,0)(1,-0.5)(1,0)
\qbezier(-0.8,-0.9)(-0.8,-0.45)(-0.2,-0.1)
\qbezier(0.2,0.1)(1,0.5)(1,0)
}\right)\  =\ 0.\
\end{equation}

Recall that one of our goals is to seek
 an intrinsic explanation for the result of
 Bar-Natan and
Garoufalidis. In this result, the invariant employed was always the
writhe-normalised quantum group invariant. If we are seeking to generalise the
scaling of the 
dimension of the chosen representation with a cabling operation, then we must
conceive of a way to remove the framing dependence once cabling has introduced
it.

What we require is an operator which projects out the subspace of chord
diagrams with isolated chords, whilst preserving 4-T relations. We will investigate the 
properties of such a projector here, \cite{swill}. First, some technical maps:

\begin{defn}
Write $s: \aaa_n 
\longrightarrow \aaa_{n-1}$, for the map which acts on chord diagrams
by summing over ways of deleting a single chord, extended linearly.
For example
\begin{equation}
s\left(
\Picture{
\FullCircle
\FullChord[10,0]
\FullChord[6,11]
\FullChord[3,9]
\FullChord[4,8]
}\right)\ =\ 
3\,\Picture{
\FullCircle
\FullChord[6,11]
\FullChord[3,9]
\FullChord[4,8]
}\ +\
\Picture{
\FullCircle
\FullChord[10,0]
\FullChord[3,9]
\FullChord[4,8]}
\ \ .
\end{equation}

\

Write $\theta:\aaa_n \longrightarrow \aaa_{n+1}$ for the map that connect-sums
in the chord diagram with a single chord.
\rtb
\end{defn}

These maps are well-defined, preserving 4-T relations. 
If we take a 4-T relation,
then we can see that the image is a sum of two terms, being terms where a chord
is removed which is either active or not in the 4-T relation. The
latter obviously still involves a 4-T relation at lesser degree, whilst the
former vanish as terms where an ``active'' chord is removed always pair up and
cancel.

It is also worth noting that with respect to the
product in the natural 
algebra, the $s$ operation 
satisfies a Leibniz rule $s(a.b)\ =\ s(a).b\ +\ a.s(b)$.
With this operation we can re-express the deframing operation in a
form that will suit our analyses. The definition follows.

\begin{defn}
Define $\phi: \aaa_n \longrightarrow \aaa_n$ to be the following operation.
\bq
\phi\ =\ Id\ -\ \theta\circ s\ +\ \frac{\theta^2\circ s^2}{2!} -\ \ldots\ +\frac{(-1)^n
\theta^n \circ s^n}{n!}.
\eq
\end{defn}

\begin{lemma}
\bq
s\circ\phi\ =\ 0.
\eq
\end{lemma}
{\bf \underline{Proof.}}

\

We have from the Leibniz rule that $s\circ\theta^m\circ s^m\ =\ m\theta^{m-1}\circ s^m\ +\
\theta^m\circ s^{m+1}.$

Consider any two successive terms in the expansion of $s\circ \phi$,
and use the fact that
\begin{eqnarray}
\ & & s\circ \frac{1}{m!}\theta^m\circ s^m\ -\
\frac{1}{(m+1)!}\theta^{m+1}\circ s^{m+1} \nonumber  \\
& = & \frac{1}{m!} m\theta^{m-1}\circ s^{m}\ + \frac{1}{m!}\theta^m\circ s^{m+1}\ -\
\frac{1}{m!}\theta^m\circ s^{m+1}\ -\ \frac{1}{(m+1)!}\theta^{m+1}\circ s^{m+2}
\nonumber \\
& = & \frac{1}{(m-1)!}\theta^{m-1}\circ s^{m}\ -\
\frac{1}{(m+1)!}\theta^{m+1}\circ s^{m+2}. 
\end{eqnarray}

Adding terms, and using the fact that
$s^{n+1}\ =\ 0$ on diagrams with $n$ chords, the lemma follows.
\rtb

\begin{corollary}

$\phi$ is a projection operator.

\end{corollary}

{\bf \underline{Proof.}}

\begin{eqnarray}
\phi^2 & = & 
(Id\ -\ \theta\circ s\ +\ \frac{\theta^2\circ s^2}{2!} -\ \ldots\ +\frac{(-1)^n
\theta^n \circ s^n}{n!})\circ \phi \nonumber \\
& = & \phi
\end{eqnarray}
\rtb

\begin{lemma}
\[
\phi \left( \Picture{
\DottedCircle
\FullChord[5,6]
\Arc[5]
\Arc[6]
}
\right)
\ \ =\ \ 0.
\]
\end{lemma}

{\bf \underline{Proof.}}

Consider
 a chord diagram with an isolated chord. $\theta^m \circ s^m$ is the 
operation
of summing over the ``trivialisation" of all choices of $m$ chords. This produces two
 terms, one where
the already trivial chord is included in the choice, and one where it is not.
 The term where it is not cancels with the term where it is included, at 
the $(m+1)$-th term in the expansion of $\phi$. This
is because they are the same diagram and there are $m+1$ more such terms from
$\theta^{m+1}\circ s^{m+1}$. This terminates at $\theta^n\circ s^n$, as there is no term
without trivial chord included in the set to be trivialised.
\rtb

Thus $\phi$ is our desired `deframing' operator.
There is a nice description of the invariant subspace of deframing over the
eigenvectors of cabling. The following criterion is useful for this task.

\begin{lemma}
For $v\, \ci\, \aaa_n,\ \phi(v) = v$ if and only if $s(v) = 0$.

\end{lemma}

{\bf \underline{Proof.}}

It is clear that if $s(v)\ =\ 0$, then $\phi(v)\ =\ v$, by 
construction. Further, we have shown that $s\circ \phi\ =\ 0$, so that if $\phi(v)\ =\
v$, then $s(v)\ =\ s(\phi(v))\ =\ 0$.
\rtb

How does the operator $s$
act on the eigenvectors of cabling? First consider
its action on Feynman diagrams. 

\begin{lemma}

The operator 
$s$ acts on FDs by summing over all ways of removing a chord. 
If there are no chords, then 
it takes the value zero.

\end{lemma}

{\bf \underline{Proof.}}

Take a Feynman diagram. Resolve all but one of the trivalent vertices.
In the sum over removals of chords, the terms where we remove the chords
resulting from the trivalent vertex cancel. Thus with one vertex, the action of $s$ is equivalent 
to a sum over removals of the non-participating chords. This understanding 
proceeds by induction.
\rtb

The above procedure
translates simply to our previously constructed cabling eigenvectors -- the operator
$s$ acts on
a SFD by striking out a single isolated chord in all possible
ways, if they 
exist, otherwise $s$ maps the SFD to zero.
For example,
\begin{equation}
s\left(
\Picture{
\put(-0.8,1){\line(0,-1){2}}
\put(-0.2,1){\line(0,-1){2}}
\put(0.8,1){\line(0,-1){2}}
\put(-0.2,0.4){\line(1,0){1}}
\put(-0.2,-0.4){\line(1,0){1}}
\put(0.3,0.4){\line(0,-1){0.8}}
}
\right)\ =\ 
\Picture{
\put(-0.2,1){\line(0,-1){2}}
\put(0.8,1){\line(0,-1){2}}
\put(-0.2,0.4){\line(1,0){1}}
\put(-0.2,-0.4){\line(1,0){1}}
\put(0.3,0.4){\line(0,-1){0.8}}
}\ ,\ \
s\left(
\Picture{
\put(-0.5,1){\line(0,-1){2}}
\put(0.5,1){\line(0,-1){2}}
}\right)\  =\ 
2\,\Picture{
\put(0,1){\line(0,-1){2}}
}\ ,\ \
s\left(
\Picture{
\put(-0.7,1){\line(0,-1){2}}
\put(0.3,1){\line(0,-1){2}}
\put(-0.7,0.4){\line(1,0){1}}
\put(-0.7,-0.4){\line(1,0){1}}
\put(-0.2,0.4){\line(0,-1){0.8}}
}\right)\ \ =\ \ 0.
\end{equation}

It is immediate that the set of SFDs with no isolated chords is
in the kernel of $s$, 
and hence in the deframing invariant subspace.
With the unique SFD at level one, these generate the full kernel of
$s$. 

Note that on the space of SFDs, the primitive subspace of the Hopf algebra 
is spanned by the connected diagrams. 

The first part of Theorem 2 now follows easily. 
Recall that the eigenvalue of a SFD
with $m$ legs 
under $\psi^n$ is $n^m$. The eigenvectors  
in the deframing invariant subspace whose eigenvalues are of leading order
are those that have 
no isolated chords, and the
maximum number of legs. 

Note first that two legs cannot join in a trivalent
vertex. This follows immediately from the
antisymmetry of trivalent vertices. (Note that the edges are not oriented, 
the arrows here point to the rest of the diagram.)

\begin{equation}
\Picture{
\put(-0.5,0){\vector(-1,0){0.25}}
\qbezier(0,0)(0.5,0.866)(0.5,0.866)
\qbezier(0,0)(0.5,-0.866)(0.5,-0.866)
\qbezier(0,0)(-0.5,0)(-0.5,0)
}\ \ =\ \ -\
\Picture{
\put(-0.5,0){\vector(-1,0){0.25}}
\qbezier(0,0)(0,-1)(0.5,0.866)
\qbezier(0,0)(0,+1)(0.5,-0.866)
\qbezier(0,0)(-0.5,0)(-0.5,0)
}\ \ =\ \ -\
\Picture{
\put(-0.5,0){\vector(-1,0){0.25}}
\qbezier(0,0)(0.5,0.866)(0.5,0.866)
\qbezier(0,0)(0.5,-0.866)(0.5,-0.866)
\qbezier(0,0)(-0.5,0)(-0.5,0)
}
\end{equation}

\

The maximum admissible number of legs then, at level $n$, is $n$, and
the eigenvectors spanning the highest weight, deframing
invariant subspace are all the ways of connecting the $n$ separate
 trivalent vertices that each leg joins, with extra edges. This corresponds
with the different ways of connecting $n$ points with closed
loops. These ways are enumerated by the partitions
of $n$.

\begin{defn}

Denote by $P$ some
partition of a positive integer $n$ (say $P=\{P_1,\ldots,P_{\# P}\}$).
We construct the SFD {\bf $\tau_{P}$}. It has $\# P$ components. The $i$th
component is a loop of $P_i$ edges with legs attached radially at every
vertex.
\rtb
\end{defn}

Some examples --

\[
\tau_{\{2\}}\ =\
\Picture{
\put(0,1){\line(0,-1){0.6}}
\put(0,-0.4){\line(0,-1){0.6}}
\qbezier(0,0.4)(-1,0)(0,-0.4)
\qbezier(0,0.4)(1,0)(0,-0.4)
}\ \ ,\ \
\tau_{\{4,2\}}\ =\
\Picture{
\put(-1,1){\line(0,-1){2}}
\put(0,1){\line(0,-1){2}}
\put(-1,0.4){\line(1,0){1}}
\put(-1,-0.4){\line(1,0){1}}
\put(1,1){\line(0,-1){0.6}}
\put(1,-0.4){\line(0,-1){0.6}}
\qbezier(1,0.4)(0,0)(1,-0.4)
\qbezier(1,0.4)(2,0)(1,-0.4)
}\ \ \ ,\ \
\tau_{\{6\}}\ =\
\Picture{
\put(-1,1){\line(0,-1){2}}
\put(1,1){\line(0,-1){2}}
\put(0,-0.4){\line(0,-1){0.6}}
\put(0,1){\line(0,-1){0.6}}
\put(-1,0.4){\line(1,0){2}}
\put(-1,-0.4){\line(1,0){2}}
}\ .\
\]

\

Observe that if such a vector is built from an odd partition 
(i.e. has an odd-legged
component) then it is zero. This again comes from the anti-symmetry condition.
We can \lq flip' such a loop over a given external chord, yielding the same 
SFD, multiplied by a minus sign - 
\[
\Picture{
\put(-1,1){\line(0,-1){0.6}}
\put(1,1){\line(0,-1){0.6}}
\put(-1,0.4){\line(1,0){2}}
\put(-1,0.4){\line(5,-4){1}}
\put(0,-0.4){\line(5,4){1}}
\put(0,-1){\line(0,1){0.6}}
}\ \ =\ \ -\,
\Picture{
\qbezier(-1,0.4)(0,0)(-1,1)
\put(1,1){\line(0,-1){0.6}}
\put(-1,0.4){\line(1,0){2}}
\put(-1,0.4){\line(5,-4){1}}
\put(0,-0.4){\line(5,4){1}}
\put(0,-1){\line(0,1){0.6}}
}\ \ =\ \
\Picture{
\qbezier(-1,0.4)(0,0)(-1,1)
\qbezier(1,0.4)(0,0)(1,1)
\put(-1,0.4){\line(1,0){2}}
\put(-1,0.4){\line(5,-4){1}}
\put(0,-0.4){\line(5,4){1}}
\put(0,-1){\line(0,1){0.6}}
}\ \ =\ \
\Picture{
\put(-1,1){\line(0,-1){0.6}}
\put(1,1){\line(0,-1){0.6}}
\put(-1,0.4){\line(1,0){2}}
\qbezier(-1,0.4)(1,0)(0,-0.4)
\qbezier(1,0.4)(-1,0)(0,-0.4)
\put(0,-1){\line(0,1){0.6}}
}\ \ =\ \ -\,
\Picture{
\put(-1,1){\line(0,-1){0.6}}
\put(1,1){\line(0,-1){0.6}}
\put(-1,0.4){\line(1,0){2}}
\put(-1,0.4){\line(5,-4){1}}
\put(0,-0.4){\line(5,4){1}}
\put(0,-1){\line(0,1){0.6}}
}
\]

\

\

It is not hard to see that moreover, the $\tau_P$ are
linearly independent (for even partitions $P$.)

We finish this section by noting that, as well as the operator $s$
defined above, there are further operations one
can perform with interesting properties. For example,
define an operator $d$ to act upon a chord diagram by summing the diagrams
obtained by replacing each chord in turn firstly by two parallel chords, and
then subtracting the diagrams obtained by replacing each chord in turn by
two intersecting chords. For example, 

\begin{equation}
d\left(
\Picture{
\FullCircle
\FullChord[3,9]
}\right)\ =\ 
\Picture{
\FullCircle
\FullChord[3,9]
\FullChord[4,8]
}\ -\
\Picture{
\FullCircle
\FullChord[3,8]
\FullChord[4,9]}\ \ ,
\end{equation}
\nonumber \\
\\
\begin{equation}
d\left(
\Picture{
\FullCircle
\FullChord[6,11]
\FullChord[3,9]
\FullChord[4,8]
}\right)\ =\ 
2\,\,\Picture{
\FullCircle
\FullChord[6,11]
\FullChord[3,9]
\FullChord[4,8]
\FullChord[2,10]
}\ -2\,\,
\Picture{
\FullCircle
\FullChord[6,11]
\FullChord[2,10]
\FullChord[3,8]
\FullChord[4,9]}
\ +\ 
\Picture{
\FullCircle
\FullChord[6,11]
\FullChord[3,9]
\FullChord[4,8]
\FullChord[7,10]
}\ -\ 
\Picture{
\FullCircle
\FullChord[7,11]
\FullChord[3,9]
\FullChord[4,8]
\FullChord[6,10]
}\ \ .
\end{equation}

\

It is straightforward (if tedious) 
to show that this operator preserves the 4-T relations. 
One can also show that the operators $d$ and $s$ map
between cabling eigenvectors, and that moreover $d\circ s - s\circ d=0$. 
There are also interesting generalisations of these operators.
Simple realisations of these operators exist for some of
the Lie algebraic weight systems. For example, in \cite{FKV} it is shown
how the Alexander-Conway weight system arises using the superalgebra
$gl(1\vert 1)$.
Before applying deframing, this weight system assigns a function of
two variables $c$ and $y$ to each chord diagram. 
One can show by induction, using the recursion relation of \cite{FKV}, that the operators
 $d$ and $s$ are realised as the differential
operators $y{\partial\over\partial c}$ and ${\partial\over\partial c}$,
respectively. Thus the action of deframing can be interpreted as the specification $c=0$.


\section{Immanents and cabling eigenvectors.}

There is another way of representing the information in a chord
diagram -- by its {\it labelled intersection graph }
(LIG) \cite{BNG}.

Consider a chord diagram $D\, \ci\, \aaa_m$. Construct a labelled graph as 
follows. The 
$m$ vertices of this graph correspond to the $m$ 
chords of the diagram, and there
is an edge connecting two vertices when the corresponding chords intersect
once in
$D$. Number the vertices according to the order of appearance going
anti-clockwise from some arbitrarily chosen point on the external loop of the
chord diagram. For example, 

\

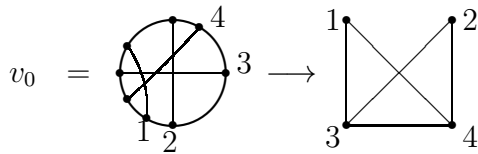
\[
v_0\ \ =\ \
\Picture{
\FullCircle
\FullChord[3,9]
\FullChord[0,6]
\FullChord[2,7]
\FullChord[5,8]
\put(-0.7,-1.3){\mbox{1}}
\put(-0.2,-1.5){\mbox{2}}
\put(1.2,0){\mbox{3}}
\put(0.7,0.866){\mbox{4}}
}\ \ \ \lra\ \
\Picture{
\put(-1,1){\line(0,-1){2}}
\put(-1,1){\line(1,-1){2}}
\put(-1,-1){\line(1,1){2}}
\put(-1,-1){\line(1,0){2}}
\put(1,1){\line(0,-1){2}}
\put(-1,1){\circle*{0.15}}
\put(-1,-1){\circle*{0.15}}
\put(1,1){\circle*{0.15}}
\put(1,-1){\circle*{0.15}}
\put(-1.4,0.8){\mbox{1}}
\put(1.2,0.8){\mbox{2}}
\put(-1.4,-1.4){\mbox{3}}
\put(1.2,-1.4){\mbox{4}}
}
\]

\

\

This information can be coded in the intersection matrix of the LIG. Construct
this matrix via the following prescription.

\begin{defn}
The {\bf Intersection Matrix} (IM) of a LIG is defined as follows.
\bq
(IM)_{ij}\ =\ \left\{ \begin{array}{ll} sign(i-j)\ & \mbox{if the vertices
labelled i and j
are 
linked}, \\ 0\ & \mbox{otherwise.} \end{array} \right.
\eq
Write $IM: \hat{\aaa}_n \rightarrow Gl(n,{\bf Z})$.
\rtb
\end{defn}
The IM for the above example follows.
\[
\left[
\begin{array}{rrrr}
0 & 0 & -1 & -1 \\
0 & 0 & -1 & -1 \\
1 & 1 & 0 & -1 \\
1 & 1 & 1 & 0 
\end{array}
\right].
\]
Up to 4-T relations, the IM contains sufficient information to reconstruct the 
original chord diagram \cite{BNG}. To build a weight system from the IM
 of the LIG,
one needs to generate numbers from the IM in a way which does not
depend
on the choice of break defining the numbering in addition to satisfying 4T
relations.
The determinant
of the IM proves a well-defined choice \cite{BNG}. The following result 
provides our connection with the AC polynomial.

\begin{fact}{\cite{BNG}}
The Alexander-Conway polynomial $C(h)$  is a series in powers of
$h$, $C(h)=\sum c_n h^n$. It is not difficult to see from the skein relation
that $c_n$ is in fact 
of finite type of order $n$. 
Choose $v\, \ci\, \aaa_n$. Then 
\begin{equation}
W_n[c_n](v)\ =\ Det(IM(v)).
\end{equation}
\rtb
\end{fact}

Immanents are alternative matrix invariants which yield well-defined
weight systems in this fashion \cite{BNG}. Denote by 
$\zsn$ the integer module generated by the
conjugacy\ classes of $S_n$ (this is not formally the group ring - we maintain the
notational conventions established in \cite{BNG}). Denote by $[\sigma]$ the conjugacy class
of $\sigma\ \ci\ S_n$.

\begin{defn}

The {\bf Universal Immanent Map} of an n$\times$n matrix $(M)_{ij}$,
 $Imm:\ Gl_n(\zzz )\ \lra\ \zsn$, is defined by 
\bq
Imm(M)\ =\ \sum_{\sigma \ci S_n} \prod_{i=1}^n (M)_{i \sigma(i)} [\sigma].
\eq
We define the universal immanent weight system $I: \aaa_n
\longrightarrow\ \zsn\ $ by $ I= Imm\circ IM$.
\rtb
\end{defn}

The conjugacy classes of $S_n$ are bijective with
the partitions of $n$. To see this, construct a graph from $\sigma \ci
S_n$ with $n$ vertices, and a link from $i$  to $\sigma(i)$. The connected
components of the resulting graph represent the corresponding partition.

To project to a ${\bf C}$-valued weight system, we compose some vector 
$W\ \ci\ \zsn^*$ (i.e. $W \ci Hom(\zsn,{\bf C})$) with $I$.
There are some distinguished elements in
$\zsn^*$. Namely, any representation of $S_n$ will furnish a well-defined
functional on conjugacy classes by taking the trace of a representative
element. With the alternating representation of
$S_n$ one obtains the usual matrix
 determinant. Taking the trivial representation,
one gets the permanent of the matrix, for example.

We can understand the universal immanent weight system differently.
First we note some graph theoretic terminology.

\begin{defn} A {\bf Hamiltonian cycle} on a graph is a directed and  non-repeating
cycle 
of at least two vertices, where consecutive vertices are linked in
the graph. 
\rtb
\end{defn}

\begin{defn} A {\bf Hamiltonian cycle decomposition} (HCD) of a graph  is a
collection of 
disjoint Hamiltonian cycles such that every vertex in the graph
appears in exactly one. 
The {\bf descent} of a HCD of a {\it labelled} graph is the number of
instances in a cycle 
decomposition where consecutive vertices in a cycle decrease in label
value.  
\rtb
\end{defn}

\begin{fact}[\cite{BNG}]
To every cycle decomposition of an $n$-verticed graph we can associate a
partition of $n$. The universal immanent invariant of a
labelled graph is precisely a sum over the partitions corresponding to the
different cycle decompositions of the graph, with each decomposition weighted
by $(-1)^d$, where d is the descent of the decomposition.
\rtb
\end{fact}
We illustrate this calculus with our
previous example -- the LIG above Definition (5.1), $v_0$. 
The cycle decompositions and descents here are
\[
\begin{array}{ll}
1\ra 4\ra 2\ra 3\ra 1 & d=2 \\
3\ra 2\ra 4\ra 1\ra 3 & d=2 \\
1\ra 4\ra 1,\ 2\ra 3\ra 2 & d=2 \\
2\ra 4\ra 2,\ 1\ra 3\ra 1 & d=2 \\
\end{array}
\]
Thus 
\[
I(v_0)\ =\ 2[4]\ +\ 2[2,2].
\]

The weighting of the decomposition by $(-1)^d$ allows us to ignore 
decompositions which include cycles of odd length: reversing the direction of
the odd-lengthed cycle produces the same decomposition with opposite
sign, which cancels in the summation.

We have an important connection between the immanent weight systems at
level $n$ and 
the highest weight deframing invariant subspace of $\aaa_n$. For each $n$, they
have precisely the same dimension, the number of possible even
partitions of $n$. This motivates the following theorem.

\newpage

\begin{theorem}
\label{theprop}
Consider $v\,\ci\, \aaa_m$.
\begin{enumerate}
\item{ If $\phi(v) =0$, then $I(v) = 0$.}
\item{ If $\phi(v)  = v$, and $\psi^n_m (v)\ =\ n^{p}v$ for $p<m$, then
\bq
I(v)=0.
\eq
}
\item{ 
\bq
I( \tau_{[\sigma]} ) = 2^{\# [\sigma]}m! [\sigma].
\eq
}
\end{enumerate}
In the above $\# [\sigma]$ denotes the number of components of $[\sigma]$.
\rtb
\end{theorem}

\

Before presenting the proof of this proposition, 
we will explain how the main
theorem follows from it. 

\begin{defn}
Take $\sigma,\rho\, \ci\, {\cal S}_n$.
Define $\delta_{[\sigma]}: \zsn \rightarrow \ccc$,
defined on the basis of $\zsn$ by
\bq
\delta_{[\sigma]}([\rho]) = \left\{ \begin{array}{l} 1\ \ \mbox{if}\
[\sigma]=[\rho], \\ 0\ \ \mbox{otherwise.} \end{array} \right.
\eq
and extend linearly. This is the canonical dual basis.

Define $\alpha_{[\sigma]}: \aaa_n \rightarrow \ccc$, by 
\bq
\alpha_{[\sigma]}=\delta_{[\sigma]}\circ I.
\eq
\rtb
\end{defn}
The $\alpha_{[\sigma]}$ span the set of immanent weight systems. Bar-Natan
and Garoufalidis showed that this subspace was equivalent to the subspace
of weight systems coming from sums of products of the coefficients of the 
AC polynomial \cite{BNG}.

\

{\bf \underline{Proof of theorem 2.}}

\

Consider the equation, for $v\, \ci\, \aaa_m$,  
\begin{eqnarray}
\widehat{\psi^{n*} W_m}(v)\ & =\ & \phi^*(\psi^{n*}
W_m)(v), \nonumber \\
& = & W_m(\psi^n_m(\phi(v))) . 
\end{eqnarray}

Expand $\phi(v)$ over the $\tau_{[\sigma]}$ and SFDs
 with
fewer univalent vertices, as $\phi(v)=\sum_{[\sigma]} b_{[\sigma]} \tau_{[\sigma]} + Rem$, where {\it Rem} denotes terms with fewer
than $m$ univalent vertices. From Theorem \ref{theprop} it then
follows that
\
\begin{eqnarray}
\alpha_{[\sigma]}(v) &=& \alpha_{[\sigma]}(\phi(v) + (Id-\phi)(v)) \\
   &=& \alpha_{[\sigma]} (\sum_{[\rho]} b_{[\rho]} \tau_{[\rho]}
      +  Rem) \\
 &=& b_{[\sigma]} 2^{\#[\sigma]}m!.
\end{eqnarray}

Applying the cabling operator and
using  Theorem 3, this implies that
 
\begin{eqnarray}
\widehat{\psi^{n*}W_m}(v) & = &
W_m\left(\psi^n_m(\sum_{[\sigma]} b_{[\sigma]} \tau_{[\sigma]}\ \ +\ \ 
 Rem) \right),
 \\
& = & n^m\sum_{[\sigma]}b_{[\sigma]} W_m(\tau_{[\sigma]})\ +\ \left(
\begin{array}{l} \mbox{lower powers} \\ \mbox{in $n$.} \end{array} \right), \\
& = & n^m\sum_{[\sigma]} k^{W_m}_{[\sigma]}
\alpha_{[\sigma]}(v)\  \ +\ \ \left( 
\begin{array}{l} \mbox{lower powers} \\ 
\mbox{in $n$.}
 \end{array} \right), 
\end{eqnarray}
setting $k^{W_m} _{[\sigma]}=1/(2^{\#[\sigma]}m!) W_m(\tau_{[\sigma]})$.
Thus 
we see  that the
$n$-th cabling of a weight system of order $m$ is a polynomial in
$n$ of highest order $n^m$, and that the coefficient of $n^m$
in this polynomial is a linear combination of immanent weight systems.
These are statements 1 and 2 of Theorem 2. 
\rtb

\
Now we consider the proof of Theorem \ref{theprop}.

\

{\bf \underline {Proof of part (1) of Theorem \ref{theprop}. }}

\

Consider $v\,\ci\, \aaa_n$
such that $\phi(v)\ =\  0$. This implies that
\bq
(Id\ -\ \theta\circ s\ +\ \frac{1}{2!}\theta^2\circ s^2 + \ldots\ +
\frac{(-1)^n}{n!}\theta^n\circ s^n )( v)\ =\ 0.
\label{above}
\eq
Now $I(w)\ =\ 0$ if $w$ has an isolated chord --  as the LIG has
 an isolated vertex, there can be no cycle decompositions. 
 Operating on both
sides of (\ref{above}) with $I$ we get
\bq
I(v)\ =\ 0.
\eq
\rtb

\section{The universal immanent weight system and FDs.}

Recall that the genus of a connected graph is calculated as 
$1-\#\mbox{vertices}+\#\mbox{edges}$. Take some FD $v$ whose
graph has $n$ components. Define $G(v)$ to be the unordered $n$-tuplet
of the genera of the connected components of the graph of $v$.
We introduce the notation $\hat{v}$ for the graph
that represents the SFD $v$. 

Our principal technical tool is:

\begin{lemma}
\label{vanish}
Take some Feynman Diagram $v\ci \aaa_n$. If the graph of $v$ has
a genus 1 component not equal to $\hat{\tau}_{\{p\}}$ for some even
integer $p$ then $I(v)$=0.
\end{lemma}

This will be proved in a later section.
Part (2) of Theorem \ref{theprop} is a statement about the values 
$I$ takes on SFDs at grade $n$ whose graphs have
less than $n$ univalent vertices. We characterise these vectors:

\begin{lemma}
Take a SFD $v\ci \aaa_n$ with less than $n$ univalent vertices. Then
$\hat{v}$ has a component of genus at least two.
\end{lemma}
This follows from a straightforward Euler characteristic calculation.
Namely: there will be at least one component of $\hat{v}$ with more 
trivalent than univalent vertices, $t>u$. That component will have genus
\begin{eqnarray*}
G & = & 1-(u+t)+(\frac{3t+u}{2}) \\
& = & 1+(\frac{t-u}{2}) > 1.
\end{eqnarray*}
 
With this understanding,
Part (2) of Theorem \ref{theprop} follows from the following lemma.

\begin{lemma}
Take a FD $v\ci \aaa_n$. If the graph of $v$ has a component of at least
genus $2$ then $I(v)=0$.
\label{highergenus}
\end{lemma}

{\bf \underline{Proof.}}

\

An STU resolution decreases the number of trivalent vertices on the
graph of $v$. $v$ is
then expressed as a linear combination of FDs whose graphs are the same.
On account of this identification it makes sense to speak
of resolving the trivalent vertices of the graph of $v$ in a particular
order.

At each resolution, the genus of the graph of the FDs in the sum can alter
in two ways. If the number of components increases by one as a result of the 
resolution then $\{g_1,g_2,\ldots,g_n\} \rightarrow \{g_1',g_1-g_1',g_2,\ldots,g_n\}$. Otherwise
(when the number of components of the graph is unchanged) 
$\{g_1,g_2,\ldots,g_n\} \rightarrow \{g_1-1,g_2,\ldots,g_n\}$.

It is always possible to resolve a choice of trivalent vertices such that $v$ 
is equal to a sum over FDs with genus $\{1,\ldots\}$. If the genus 1 component
is not $\hat{\tau}_{\{p\}}$ for some even $p$, then $I(v)=0$ from 
Lemma \ref{vanish}.

Assume then, that we have expressed (by some sequence of STU resolutions
of $v$) $v$ as a sum over FDs whose graphs have a genus 1 component 
$\tau_{\{p\}}$ for some even $p$.
The step which led to this was either $\{g_1,g_2,\ldots,g_n\}\rightarrow \{1,g_1-1,g_2,\ldots,g_n\}$ or 
$\{2,g_2,\ldots,g_n\}\rightarrow \{1,g_2,\ldots,g_n\}$. 
We show here that in both these cases we can always choose a {\it different}
sequence 
of vertex resolutions so that $v$ is expressible as a sum of FDs with 
genus 1 components {\it not} some $\hat{\tau}_{\{p\}}$ (and hence $I(v)$ vanishes
by Lemma (\ref{vanish})).

Take the first case then, where some genus $g>1$ component ``splits'' into
a genus 1 and a genus $g-1$ component when some joining vertex is
resolved. For example:

\begin{eqnarray*}
a & = & 
\Picture{
\FullCircle
\FullChord[4,9]
\InOut[0,0]
\InIn[0,1]
\InIn[0,4]
\InOut[2,5]
\InOut[3,7]
\qbezier(0.206,0.634)(-0.1,0.8)(-0.539,0.392)
\qbezier(0.206,0.634)(-0.1,0.2)(-0.539,0.392)
\qbezier(0.206,-0.634)(-0.1,-0.2)(-0.539,-0.391)
\qbezier(0.206,-0.634)(-0.1,-0.8)(-0.539,-0.391)
\Endpoint[0]\Endpoint[5]\Endpoint[7]
} \\
& & \\
& & \\
& = & 
\Picture{
\FullCircle
\FullChord[4,9]
\InOut[1,1]
\InOut[4,11]
\InOut[2,5]
\InOut[3,7]
\qbezier(0.206,0.634)(-0.1,0.8)(-0.539,0.392)
\qbezier(0.206,0.634)(-0.1,0.2)(-0.539,0.392)
\qbezier(0.206,-0.634)(-0.1,-0.2)(-0.539,-0.391)
\qbezier(0.206,-0.634)(-0.1,-0.8)(-0.539,-0.391)
\Endpoint[1]\Endpoint[11]\Endpoint[5]\Endpoint[7]
}\ \ -\ \
\Picture{
\FullCircle
\FullChord[4,9]
\InOut[1,11]
\InOut[4,1]
\InOut[2,5]
\InOut[3,7]
\qbezier(0.206,0.634)(-0.1,0.8)(-0.539,0.392)
\qbezier(0.206,0.634)(-0.1,0.2)(-0.539,0.392)
\qbezier(0.206,-0.634)(-0.1,-0.2)(-0.539,-0.391)
\qbezier(0.206,-0.634)(-0.1,-0.8)(-0.539,-0.391)
\Endpoint[1]\Endpoint[11]\Endpoint[5]\Endpoint[7]
} \\
& & \\
& \equiv & b\ \ -\ \ c.\\
G(a) & = & \{ 0, 2\} \\
G(b)=G(c) & = & \{ 0,1,1 \}.
\end{eqnarray*}

We can always choose to resolve all the vertices that make up
the genus $g-1$ subgraph instead of the ``joining'' vertex. 
The genus 1 component is always then $\hat{\tau}_{\{p\}}$ with some 
tree adjoined. Consider our example:

\[
a\ \  =\ \  
\Picture{
\FullCircle
\FullChord[4,9]
\InOut[0,0]
\InIn[0,1]
\InIn[0,4]
\InOut[2,5]
\InOut[3,7]
\qbezier(0.206,0.634)(-0.1,0.8)(-0.539,0.392)
\qbezier(0.206,0.634)(-0.1,0.2)(-0.539,0.392)
\qbezier(0.206,-0.634)(-0.1,-0.2)(-0.539,-0.391)
\qbezier(0.206,-0.634)(-0.1,-0.8)(-0.539,-0.391)
\Endpoint[0]\Endpoint[5]\Endpoint[7]
}\ \  =\ \   
 2 \,\Picture{
\FullCircle
\FullChord[4,9]
\InOut[0,0]
\InIn[0,1]
\InIn[0,4]
\InOut[1,2]
\InOut[1,3]
\InOut[3,7]
\qbezier(0.206,-0.634)(-0.1,-0.2)(-0.539,-0.391)
\qbezier(0.206,-0.634)(-0.1,-0.8)(-0.539,-0.391)
\Endpoint[0]\Endpoint[2]\Endpoint[7]\Endpoint[3]
} 
\]

\

The second possibility is that the genus 1 component $\hat{\tau}_{\{p\}}$ is 
obtained by resolving the vertex on some genus 2 component.
For example:

\begin{eqnarray*}
a & = &
\Picture{
\FullCircle
\qbezier(1,0)(1,0)(0.667,0)
\InIn[0,1]
\InIn[0,4]
\InIn[1,4]
\InIn[1,2]
\InIn[4,3]
\InIn[2,3]
\InOut[2,5]
\InOut[3,7]
\FullChord[9,3]
\FullChord[2,10]
\Endpoint[5]
\Endpoint[7]
\Endpoint[0]
} \ \\
& & \\
& & \\
& = &
\Picture{
\Endpoint[5]\Endpoint[7]\Endpoint[1]\Endpoint[11]
\FullCircle
\InOut[4,11]
\InOut[1,1]
\InIn[1,4]
\InIn[1,2]
\InIn[4,3]
\InIn[2,3]
\InOut[2,5]
\InOut[3,7]
\FullChord[9,3]
\FullChord[2,10]
}\ \ -\ \
\Picture{
\Endpoint[5]\Endpoint[7]\Endpoint[1]\Endpoint[11]
\FullCircle
\InOut[4,1]
\InOut[1,11]
\InIn[1,4]
\InIn[1,2]
\InIn[4,3]
\InIn[2,3]
\InOut[2,5]
\InOut[3,7]
\FullChord[9,3]
\FullChord[2,10]
} \\
& & \\
& & \\
& \equiv & b\ \ -\ \ c, \\
G(a) & = & \{0,0,2 \} \\
G(b)=G(c) & = & \{0,0,1 \}.
\end{eqnarray*}

Possible genus 2 components of this sort are $\hat{\tau}_{\{p\}}$ with
two of the legs joined in an extra trivalent vertex. The case $p=2$ is 
resolved 
from a null vector:

\[
\Picture{
\FullCircle
\InOut[0,0]
\qbezier(0.667,0)(0.2,0.5)(0,0.5)
\qbezier(0.667,0)(0.2,-0.5)(0,-0.5)
\put(0,-0.5){\line(0,1){1}}
\qbezier(0,0.5)(-0.2,0.5)(-0.667,0)
\qbezier(0,-0.5)(-0.2,-0.5)(-0.667,0)
}\ \ =\ \ -\
\Picture{
\FullCircle
\InOut[0,0]
\qbezier(0.667,0)(0.5,-0.5)(0,0.5)
\qbezier(0.667,0)(0.5,0.5)(0,-0.5)
\put(0,-0.5){\line(0,1){1}}
\qbezier(0,0.5)(-0.2,0.5)(-0.667,0)
\qbezier(0,-0.5)(-0.2,-0.5)(-0.667,0)
}\ \ =\ \ -\
\Picture{
\FullCircle
\InOut[0,0]
\qbezier(0.667,0)(0.2,0.5)(0,0.5)
\qbezier(0.667,0)(0.2,-0.5)(0,-0.5)
\put(0,-0.5){\line(0,1){1}}
\qbezier(0,0.5)(-0.2,0.5)(-0.667,0)
\qbezier(0,-0.5)(-0.2,-0.5)(-0.667,0)
}.
\]

\
 
Thus we can assume $p\geq 4$. Joining the two legs in this fashion partitions
the remaining legs around the loop into two sets, according to how they
appear on the internal loop. If one set has $q$ legs, the other will have
$p-q-2$ legs. If either of these sets has more than one leg, then 
choosing instead to resolve a trivalent vertex along some leg from that set 
yields a genus 1 component not some $\hat{\tau}_{\{p\}}$. Taking our example:

\begin{eqnarray*}
 &  &
\Picture{
\put(-1.4,0.4){\mbox{x}}
\put(-1.4,-0.8){\mbox{y}}
\FullCircle
\qbezier(1,0)(1,0)(0.667,0)
\InIn[0,1]
\InIn[0,4]
\InIn[1,4]
\InIn[1,2]
\InIn[4,3]
\InIn[2,3]
\InOut[2,5]
\InOut[3,7]
\FullChord[9,3]
\FullChord[2,10]
\Endpoint[5]
\Endpoint[7]
\Endpoint[0]
} \ \\
& & \\
& & \\
& = &
\Picture{
\Endpoint[6]
\FullCircle
\qbezier(1,0)(1,0)(0.667,0)
\InIn[0,1]
\InIn[0,4]
\InIn[1,4]
\InIn[4,3]
\InOut[1,5]
\InOut[3,6]
\InOut[3,7]
\FullChord[9,3]
\FullChord[2,10]
\Endpoint[5]
\Endpoint[7]
\Endpoint[0]
} \ \ -\ \
\Picture{
\Endpoint[6]
\FullCircle
\qbezier(1,0)(1,0)(0.667,0)
\InIn[0,1]
\InIn[0,4]
\InIn[1,4]
\InIn[4,3]
\InOut[1,6]
\InOut[3,5]
\InOut[3,7]
\FullChord[9,3]
\FullChord[2,10]
\Endpoint[5]
\Endpoint[7]
\Endpoint[0]
}.
\end{eqnarray*}

\

In this case we see that the extra join divides the other legs of $\hat{\tau}_{\{4\}}$ into sets $\{x,y\}$ and $\phi$ (the joined legs are adjacent on the loop).
As $\{x,y\}$ has two legs, resolving along
the leg $x$ yields $\hat{\tau}_{3}$ with a tree adjoined.

There is one exceptional case. Joining opposite legs of $\hat{\tau}_{\{4\}}$ in an
extra vertex gives two sets of one leg each (resolving along any vertex always
leads to $\tau_{\{4\}}$ which we cannot say $I$ vanishes on). Happily, this
possibility is the zero vector.

\

\begin{eqnarray}
\Picture{
\FullCircle
\Endpoint[3]\Endpoint[5]\Endpoint[10]
\InOut[1,3]
\InOut[2,5]
\InOut[4,10]
\InIn[0,1]
\InIn[0,2]
\InIn[0,4]
\InIn[3,4]
\InIn[1,3]
\InIn[2,3]
} & = & 
\ -\ \Picture{
\FullCircle
\Endpoint[3]\Endpoint[5]\Endpoint[10]
\InOut[1,3]
\InOut[2,5]
\InOut[4,10]
\InIn[0,2]
\InIn[0,4]
\InIn[3,4]
\InIn[2,3]
\qbezier(0.206,0.634)(-0.1,0.634)(0.667,0)
\qbezier(-0.539,-0.391)(0.5,0.634)(0.206,0.634)
}\ \ =\ \
\Picture{
\FullCircle
\Endpoint[3]\Endpoint[5]\Endpoint[10]
\InOut[1,3]
\InOut[2,5]
\InOut[4,10]
\InIn[0,2]
\InIn[2,3]
\qbezier(0.206,0.634)(-0.1,0.634)(0.667,0)
\qbezier(-0.539,-0.391)(0.5,0.634)(0.206,0.634)
\qbezier(0.206,-0.634)(-0.3,-0.634)(0.667,0)
\qbezier(0.206,-0.634)(0.4,-0.4)(-0.539,-0.391)
} \nonumber \\
& & \nonumber \\
& & \nonumber \\
& \equiv & \ \
\Picture{
\FullCircle
\Endpoint[3]\Endpoint[5]\Endpoint[10]
\InOut[1,3]
\InOut[2,5]
\InOut[4,10]
\InIn[0,1]
\InIn[0,4]
\InIn[3,4]
\InIn[1,3]
\qbezier(-0.539,0.392)(-0.7,0.15)(0.667,0)
\qbezier(-0.539,-0.391)(-0.2,0.392)(-0.539,0.392)
}\ \ =\ -\
\Picture{
\FullCircle
\Endpoint[3]\Endpoint[5]\Endpoint[10]
\InOut[1,3]
\InOut[2,5]
\InOut[4,10]
\InIn[0,1]
\InIn[0,2]
\InIn[0,4]
\InIn[3,4]
\InIn[1,3]
\InIn[2,3]
}. \nonumber \\
& & 
\end{eqnarray}
\rtb

We turn to the proof of Part 3 of Theorem \ref{theprop}.
Write

\begin{eqnarray}
\setlength{\unitlength}{40pt}
v & = & 
\Picture{
\put(-0.6,0.3){\mbox{a}}
\put(0.2,0.3){\mbox{b}}
\put(-1,-0.5){\line(1,0){0.4}}
\put(-0.2,-0.5){\line(1,0){0.4}}
\put(0.6,-0.5){\line(1,0){0.2}}
\put(1.0,-0.5){\vector(1,0){0.4}}
\qbezier[4](-0.6,-0.5)(-0.6,-0.5)(-0.2,-0.5)
\qbezier[4](0.2,-0.5)(0.2,-0.5)(0.6,-0.5)
\put(-0.15,-0.6){\framebox(0.3,0.2)}
\put(-0.4,-0.5){\oval(0.6,1.0)[t]}
\put(0.4,-0.5){\oval(0.6,1.0)[t]}
}
\nonumber \\
& & \nonumber \\
& = &
\Picture{
\put(-1,-0.5){\line(1,0){0.4}}
\put(-0.2,-0.5){\line(1,0){0.4}}
\put(0.6,-0.5){\line(1,0){0.2}}
\put(1.0,-0.5){\vector(1,0){0.4}}
\qbezier[4](-0.6,-0.5)(-0.6,-0.5)(-0.2,-0.5)
\qbezier[4](0.2,-0.5)(0.2,-0.5)(0.6,-0.5)
\put(-0.4,-0.5){\oval(0.6,1.0)[t]}
\put(0.4,-0.5){\oval(0.6,1.0)[t]}
}\ \ \ \
-\ \ \ \
\Picture{
\put(-1,-0.5){\line(1,0){0.4}}
\put(-0.2,-0.5){\line(1,0){0.4}}
\put(0.6,-0.5){\line(1,0){0.2}}
\put(1.0,-0.5){\vector(1,0){0.4}}
\qbezier[4](-0.6,-0.5)(-0.6,-0.5)(-0.2,-0.5)
\qbezier[4](0.2,-0.5)(0.2,-0.5)(0.6,-0.5)
\put(-0.3,-0.5){\oval(0.8,1.0)[t]}
\put(0.3,-0.5){\oval(0.8,1.0)[t]}
}\nonumber \\
& & \nonumber \\
& = & v_L - v_R
\setlength{\unitlength}{20pt}
\end{eqnarray}

\begin{lemma}
The cycle decomposition sum of the labelled intersection graph of
$v$ (above) is a linear combination of cycle decompositions, each of
which includes
a step from the vertex corresponding to $a$ to the vertex corresponding to $b$.
\label{muststep}
\end{lemma}

{\bf \underline{Proof.}}

\

On the level of labelled intersection graphs, $v_L-v_R$ looks like:

\begin{equation}
\Picture{
\put(-0.5,0){\circle*{0.15}}
\put(0.5,0){\circle*{0.15}}
\put(-0.5,0){\line(-1,-1){0.5}}
\put(-0.5,0){\line(-1,1){0.5}}
\put(-0.5,0){\line(-1,0){0.75}}
\put(0.5,0){\line(1,1){0.5}}
\put(0.5,0){\line(1,-1){0.5}}
}\ \ \ \ -\ \ \ \
\Picture{
\put(-0.5,0){\circle*{0.15}}
\put(0.5,0){\circle*{0.15}}
\put(-0.5,0){\line(-1,-1){0.5}}
\put(-0.5,0){\line(-1,1){0.5}}
\put(-0.5,0){\line(-1,0){0.75}}
\put(0.5,0){\line(1,1){0.5}}
\put(0.5,0){\line(1,-1){0.5}}
\put(-0.5,0){\line(1,0){1}}
}
\end{equation}

\

(where we have drawn in possible further links).
Note that the labels on the vertices corresponding to $a$ and $b$ are
either unchanged between $v_L$ and $v_R$,
or they swap (in which case the labels are consecutive).  
Any cycle decomposition of $v_R$ which does not include a step from the
vertex corresponding to chord $a$ to  the vertex corresponding to chord $b$
decomposes $v_L$ with the same descent, and hence cancels.
\rtb

\

Define $\tau_{P}'$ to be the FD that corresponds to the planar
embedding of the graph $\hat{\tau}_{P}$ in a Wilson loop. For example:

\begin{equation}
\tau_{\{4\}}'\  = \ 
\Picture{
\FullCircle
\InIn[1,2]
\InIn[2,3]
\InIn[3,4]
\InIn[4,1]
\InOut[1,2]
\InOut[2,5]
\InOut[3,7]
\InOut[4,10]
}.
\end{equation}

\

\begin{lemma}
\bq
I(\tau_{\{p\}}') = 2 \{ p \}.
\eq
\label{basicloop}
\end{lemma}
{\bf \underline{Proof.}}
Take the example $\tau_{\{4\}}'$. Resolve every trivalent vertex through the
leg it joins:
\begin{equation}
\setlength{\unitlength}{30pt}
\tau_{\{4\}}' = 
\ \ \ \ \ \ \Picture{
\put(-2.0,-0.5){\vector(1,0){3.2}}
\qbezier[8](-1.4,-0.5)(-1.4,-0.5)(-1.0,-0.5)
\qbezier[8](-0.6,-0.5)(-0.6,-0.5)(-0.2,-0.5)
\qbezier[8](0.2,-0.5)(0.2,-0.5)(0.6,-0.5)
\qbezier[8](1.0,-0.5)(1.0,-0.5)(1.4,-0.5)
\put(-1.75,-0.6){\framebox(0.3,0.2)}
\put(-0.95,-0.6){\framebox(0.3,0.2)}
\put(-0.15,-0.6){\framebox(0.3,0.2)}
\put(0.65,-0.6){\framebox(0.3,0.2)}
\put(-1.2,-0.5){\oval(0.6,1.0)[t]}
\put(-0.4,-0.5){\oval(0.6,1.0)[t]}
\put(0.4,-0.5){\oval(0.6,1.0)[t]}
\put(-0.4,-0.5){\oval(2.6,2.0)[t]}
}
\setlength{\unitlength}{20pt}
\end{equation}

\

\

Lemma \ref{muststep} indicates that the cycle 
decomposition sum is a linear combination of cycle decompositions which
cycle around the four chords in the order in which they
meet. The cycle decompositions are:

\begin{equation}
\Picture{
\put(-1.2,-0.5){\vector(1,0){2.4}}
\put(0,-0.5){\oval(1.6,2.4)[t]}
\put(-0.65,-0.5){\oval(0.7,1.2)[t]}
\put(0.65,-0.5){\oval(0.7,1.2)[t]}
\put(0,-0.5){\oval(1,1.2)[t]}
\put(-1.2,-1){\mbox{1}}
\put(-1,-1){\mbox{2}}
\put(-0.7,-1){\mbox{3}}
\put(0.1,-1){\mbox{4}}
}\ \ \ \ ,\ \ \ \
\begin{array}{l}
1\rightarrow 3\rightarrow 4\rightarrow 2\rightarrow 1, \\
1\rightarrow 2\rightarrow 4\rightarrow 3\rightarrow 1.
\end{array}
\end{equation}

\

Both these cycles have descent two.
Thus $I(\tau_{\{4\}}')=2\{4\}$. The only case that is treated slightly
differently to this example is for $\tau_{\{2\}}'$. Here:

\begin{equation}
\tau_{\{2\}}' = 2\left( \ \
\Picture{
\put(-1,-0.5){\line(1,0){0.4}}
\put(-0.2,-0.5){\line(1,0){0.4}}
\put(-1,-0.5){\vector(1,0){2.4}}
\qbezier[4](-0.6,-0.5)(-0.6,-0.5)(-0.2,-0.5)
\qbezier[4](0.2,-0.5)(0.2,-0.5)(0.6,-0.5)
\put(-0.4,-0.5){\oval(0.6,1.0)[t]}
\put(0.4,-0.5){\oval(0.6,1.0)[t]}
}\ \ \ \
-\ \ \ \
\Picture{
\put(-1,-0.5){\line(1,0){0.4}}
\put(-0.2,-0.5){\line(1,0){0.4}}
\put(-1,-0.5){\vector(1,0){2.4}}
\qbezier[4](-0.6,-0.5)(-0.6,-0.5)(-0.2,-0.5)
\qbezier[4](0.2,-0.5)(0.2,-0.5)(0.6,-0.5)
\put(-0.3,-0.5){\oval(0.8,1.0)[t]}
\put(0.3,-0.5){\oval(0.8,1.0)[t]}
}\ \ \ \ \right).
\end{equation}

\

It is easy to see in this case that $I(\tau_{\{2\}}')=2\{2\}.$
\rtb

\

Recall that $\zsn$ has a generator for every partition of $n$. We define
a multiplication 
$\cdot:{\cal Z}S_n \times {\cal Z}S_m \rightarrow {\cal Z}S_{n+m}$
defined on the generators by juxtaposition of partitions (i.e. $\{p_1,\ldots,p_i\}\cdot \{q_1,\ldots,q_j\} = \{p_1,\ldots,p_i,q_1,\ldots,q_j\}$), extended
linearly. The following property is manifest from the definition of $I$.
\begin{lemma}
Take $v\ci \aaa_n,\ w\ci \aaa_m$. Then
\bq
I(v.w) = I(v)\cdot I(w).
\eq
\label{mult}
\end{lemma}

Take some SFD $\tau_{\{p_1,\ldots,p_i\}}$. Recall that this is a
sum over the FDs that correspond to all different orderings of the legs
of the graph $\hat{\tau}_{\{p_1,\ldots,p_i\}}$ on a Wilson loop. Recall
that on account of the STU relations, the FDs corresponding to different 
orderings of the legs on the Wilson loop differ by FDs with more internal 
vertices. In fact here we observe:
\bq
\tau_{\{p_1,\ldots,p_i\}} = (p_1+\ldots+p_i)!\tau_{\{p_1\}}' . \ldots . \tau_{\{p_i\}}' +
\left\{ \begin{array}{l}
\mbox{FDs of genus} \\
\geq 2
\end{array}
\right\}.
\eq

This observation, together with Lemma \ref{highergenus},  Lemma \ref{basicloop}
and Lemma \ref{mult} yields Part 3 of Theorem \ref{theprop}.

\section{Proof of Lemma \ref{vanish}}

\

Lemma \ref{vanish} details conditions under which $I$ vanishes. We separate
the proof into two parts.

\

{\bf {Lemma \ref{vanish} A}}

\

{\it If the graph of a FD $v$ has a component of its graph $\hat{\tau}_{\{p\}}$ for some odd integer $p$, then $I(v)=0$. }

\

{\bf \underline{Proof.}}

\

Such a FD has a presentation, taking the example $p=3$,

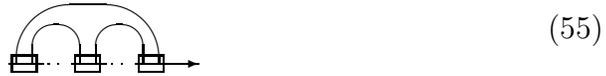
\begin{equation}
\setlength{\unitlength}{30pt}
\Picture{
\put(-1,-0.5){\line(1,0){0.4}}
\put(-0.2,-0.5){\line(1,0){0.4}}
\put(0.6,-0.5){\vector(1,0){0.8}}
\qbezier[4](-0.6,-0.5)(-0.6,-0.5)(-0.2,-0.5)
\qbezier[4](0.2,-0.5)(0.2,-0.5)(0.6,-0.5)
\put(-0.95,-0.6){\framebox(0.3,0.2)}
\put(-0.15,-0.6){\framebox(0.3,0.2)}
\put(0.65,-0.6){\framebox(0.3,0.2)}
\put(-0.4,-0.5){\oval(0.6,1.0)[t]}
\put(0.4,-0.5){\oval(0.6,1.0)[t]}
\put(0,-0.5){\oval(1.8,1.5)[t]}
}
\setlength{\unitlength}{20pt}
\end{equation}

\

By Lemma \ref{muststep} the cycle decomposition sum is a linear combination of 
cycle decompositions which include a 3-cycle around these chords. However, all
decompositions with odd cycles cancel on account of the weighting by descent.
\rtb

\

{\bf Lemma \ref{vanish} B.}

\

{\it If the graph of a FD $v$ has a genus 1 component not some 
$\hat{\tau}_{\{p\}}$ then $I(v)=0$.}

\

{\bf \underline{Proof.}}

\

A connected genus one trivalent graph $v$ has a single cycle with a 
number of 'trees'
attached (with the obvious meaning):

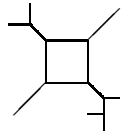
\[
\Picture{
\put(-0.4,-0.4){\line(1,0){0.8}}
\put(-0.4,-0.4){\line(0,1){0.8}}
\put(-0.4,0.4){\line(1,0){0.8}}
\put(0.4,-0.4){\line(0,1){0.8}}
\put(-0.4,-0.4){\line(-1,-1){0.6}}
\qbezier(-0.7,0.7)(-0.7,0.7)(-0.4,0.4)
\put(-0.7,0.7){\line(-1,0){0.4}}
\put(-0.7,0.7){\line(0,1){0.4}}
\put(0.4,0.4){\line(1,1){0.6}}
\qbezier(0.4,-0.4)(0.4,-0.4)(0.7,-0.7)
\put(0.7,-0.7){\line(1,0){0.3}}
\put(0.7,-0.7){\line(0,-1){0.6}}
\put(0.7,-1){\line(-1,0){0.3}}
}
\]

\

Via a sequence of STU resolutions any FD $v$ containing a genus 1
component not some $\hat{\tau}_{\{p\}}$ may be expressed as a linear
combination of FDs whose graphs contain a component from the 
following list (i.e. some $\hat{\tau}_{\{p\}}$ with one extra trivalent 
vertex):

\

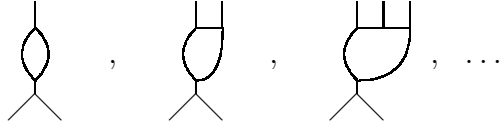
\[
\Picture{
\put(0,-0.5){\line(0,-1){0.25}}
\put(0,-0.75){\line(1,-1){0.5}}
\put(0,-0.75){\line(-1,-1){0.5}}
\put(0,0.5){\line(0,1){0.5}}
\qbezier(0,0.5)(0.5,0)(0,-0.5)
\qbezier(0,0.5)(-0.5,0)(0,-0.5)
}
\ \ ,\ \ 
\Picture{
\put(0,-0.5){\line(0,-1){0.25}}
\put(0,-0.75){\line(1,-1){0.5}}
\put(0,-0.75){\line(-1,-1){0.5}}
\put(0,0.5){\line(0,1){0.5}}
\qbezier(0,0.5)(-0.5,0)(0,-0.5)
\put(0,0.5){\line(1,0){0.5}}
\put(0.5,0.5){\line(0,1){0.5}}
\qbezier(0,-0.5)(0.5,-0.5)(0.5,0.5)
}
\ \ ,\ \
\Picture{
\put(0,-0.5){\line(0,-1){0.25}}
\put(0,-0.75){\line(1,-1){0.5}}
\put(0,-0.75){\line(-1,-1){0.5}}
\put(0,0.5){\line(0,1){0.5}}
\qbezier(0,0.5)(-0.5,0)(0,-0.5)
\put(0,0.5){\line(1,0){0.5}}
\put(0.5,0.5){\line(0,1){0.5}}
\put(0.5,0.5){\line(1,0){0.5}}
\put(1,0.5){\line(0,1){0.5}}
\qbezier(0,-0.5)(1,-0.5)(1,0.5)
}
\ \ ,\ \ \ldots 
\]

\

\

Thus we need to show that $I(v)=0$ when $v$ has such a component in its graph.
There is no loss of generality if we assume there are no trivalent vertices 
in the other graph components (i.e. they are all chords). There are two
cases to address: when the cycle has two edges, and when the cycle has more 
than two edges.

Assume first that the cycle has more than two edges. Such a FD has a
presentation as follows:

\begin{eqnarray}
\Picture{
\DottedCircle
\Endpoint[1]
\Endpoint[11]
\Endpoint[5]
\Endpoint[7]
\Endpoint[9]
\Arc[1]
\Arc[11]
\Arc[5]
\Arc[7]
\Arc[9]
\InOut[0,1]
\InOut[0,11]
\InIn[0,5]
\InIn[5,2]
\InIn[2,3]
\InIn[3,4]
\InIn[4,5]
\InOut[2,5]
\InOut[3,7]
\InOut[4,9]}
 & \rightarrow &
\setlength{\unitlength}{30pt}
\ \ \ \ \ \ \ \ \Picture{
\put(-1.75,0.5){\mbox{a}}
\put(-1.3,0.1){\mbox{b}}
\put(-0.5,0.1){\mbox{c}}
\put(0.3,0.1){\mbox{d}}
\put(1.2,-0.2){\mbox{*}}
\put(-2.0,-0.5){\line(1,0){0.6}}
\put(-1.0,-0.5){\line(1,0){0.4}}
\put(-0.2,-0.5){\line(1,0){0.4}}
\put(0.6,-0.5){\line(1,0){0.62}}
\put(1.4,-0.5){\vector(1,0){0.4}}
\qbezier[8](-1.4,-0.5)(-1.4,-0.5)(-1.0,-0.5)
\qbezier[8](-0.6,-0.5)(-0.6,-0.5)(-0.2,-0.5)
\qbezier[8](0.2,-0.5)(0.2,-0.5)(0.6,-0.5)
\qbezier[8](1.0,-0.5)(1.0,-0.5)(1.4,-0.5)
\put(-1.75,-0.6){\framebox(0.3,0.2)}
\put(-0.95,-0.6){\framebox(0.3,0.2)}
\put(-0.15,-0.6){\framebox(0.3,0.2)}
\put(0.65,-0.6){\framebox(0.3,0.2)}
\put(-1.2,-0.5){\oval(0.6,1.0)[t]}
\put(-0.4,-0.5){\oval(0.6,1.0)[t]}
\put(0.4,-0.5){\oval(0.6,1.0)[t]}
\put(-0.4,-0.5){\oval(2.6,2.0)[t]}
\put(1.3,-0.5){\oval(0.4,0.5)[t]}
}\ \ \ \ \ \ \ \ -\ \ \ \ \ \ \ \ \Picture{
\put(-2.0,-0.5){\line(1,0){0.6}}
\put(-1.0,-0.5){\line(1,0){0.4}}
\put(-0.2,-0.5){\line(1,0){0.4}}
\put(0.38,-0.5){\line(1,0){0.62}}
\put(1.4,-0.5){\vector(1,0){0.4}}
\qbezier[8](-1.4,-0.5)(-1.4,-0.5)(-1.0,-0.5)
\qbezier[8](-0.6,-0.5)(-0.6,-0.5)(-0.2,-0.5)
\qbezier[8](0.2,-0.5)(0.2,-0.5)(0.6,-0.5)
\qbezier[8](1.0,-0.5)(1.0,-0.5)(1.4,-0.5)
\put(-1.75,-0.6){\framebox(0.3,0.2)}
\put(-0.95,-0.6){\framebox(0.3,0.2)}
\put(-0.15,-0.6){\framebox(0.3,0.2)}
\put(0.65,-0.6){\framebox(0.3,0.2)}
\put(-1.2,-0.5){\oval(0.6,1.0)[t]}
\put(-0.4,-0.5){\oval(0.6,1.0)[t]}
\put(0.4,-0.5){\oval(0.6,1.0)[t]}
\put(-0.4,-0.5){\oval(2.6,2.0)[t]}
\put(1,-0.5){\oval(1,0.5)[t]}
}\ \ \ \, \nonumber \\
& & \nonumber \\
& & \nonumber \\
& = & \ \ v_L \ \ -\ \ v_R.
\setlength{\unitlength}{30pt}
\end{eqnarray}

\

From Lemma \ref{muststep} we see that the cycle decomposition sum of
$v_L$ (or $v_R$)  is a linear combination  of cycle decompositions which
include a cycle around the vertices corresponding to the chords a,b,c and d.

However note further, that any cycle decomposition which does not include a 
step from the vertex corresponding to * to at least one of a or d appears
equally signed in the cycle decomposition sums of $v_L$ and $v_R$ and hence
cancels in the sum (one must be careful to check that the descents are the
same, the point being that when labels swap, they are consecutive). 

Thus, the cycle decomposition sum is zero. $I(v)=0$.

The logic for when the genus 1 component has a two cycle is almost
identical. In this case $v$ has a presentation:

\begin{eqnarray}
\Picture{
\DottedCircle
\qbezier(-1,0)(-1,0)(-0.4,0)
\qbezier(-0.4,0)(0,0.5)(0.4,0)
\qbezier(-0.4,0)(0,-0.5)(0.4,0)
\qbezier(1,0)(1,0)(0.4,0)
\InOut[0,2]
\Endpoint[0]\Arc[0]
\Endpoint[6]\Arc[6]
\Endpoint[2]\Arc[2]
} & \rightarrow
&
\setlength{\unitlength}{30pt}
\ \ \ \ 
2\left\{
\Picture{
\put(-1,-0.5){\line(1,0){0.8}}
\put(0.2,-0.5){\line(1,0){0.4}}
\put(1,-0.5){\vector(1,0){0.4}}
\qbezier[8](-0.2,-0.5)(-0.2,-0.5)(0.2,-0.5)
\qbezier[8](0.6,-0.5)(0.6,-0.5)(1,-0.5)
\put(-0.55,-0.6){\framebox(0.3,0.2)}
\put(0.25,-0.6){\framebox(0.3,0.2)}
\put(-0.65,-0.5){\oval(0.3,0.6)[t]}
\put(0,-0.5){\oval(0.6,1.0)[t]}
\put(0.85,-0.5){\oval(0.7,1.0)[t]}}
\ \ \ \ \right\}
\setlength{\unitlength}{20pt}
\end{eqnarray}
\rtb

\section{Conclusions}
 
The result we have proved in Theorem 2
 is a statement about weight systems. What this says
about the actual knot invariants is more subtle. Recall that the
Kontsevich integral inverts weight systems \cite{Kon,BN1}.
 That is, it is a knot invariant
taking values in the algebra of chord diagrams $Z_K: \{ knots \} \lra \aaa_*.$
In general we know that for a knot invariant of finite order $m$ 
\begin{equation}
W_m(V)\circ Z_K\ =\ V\ +\ \left(
\begin{array}{l}
\mbox{invariants} \\
\mbox{of order $<m$}
\end{array}
\right).
\label{equivthing}
\end{equation}

This expression can be iterated to construct a weight system: 
\begin{eqnarray}
\hat{W}^V &  = & \hat{W}^V_m\ +\
\hat{W}^V_{m-1}\ +\ \ldots\ +\ \hat{W}^V_{1}\ \neq\ W[V]\nonumber \\
 V & = & \hat{W}^V\circ Z_K.
\end{eqnarray} 

Consider the invariant $V^n\ =\ \hat{W}^V\circ\psi^n\circ\phi\circ Z_K$.
Obviously $V^1=V$. In this work we have shown that $V^n$ is a finite
polynomial in $n$ of highest order $n^m$ and that the coefficient of $n^m$ is
a linear sum of the knot invariants $\alpha_{[\sigma]}\circ Z_K$. The weight
systems $\alpha_{[\sigma]}$ are in the algebra of the (normalised)
Alexander-Conway weight systems \cite{BNG}. As
 these weight systems are canonical (the remainder in eqn. (\ref{equivthing})
vanishes) we have shown
that the highest order term of an arbitrary Vassiliev
knot invariant is in the algebra of coefficients of the Conway polynomial.  

This is an intruiging result: every Vassiliev knot invariant has
a term which can be calculated from traditional methods of algebraic topology.
An immediate question is whether such an understanding extends to the lower
powers in $n$ of an invariant. Alternatively, how must the Alexander construction
be perturbed to account for the next-to-highest orders?

There are some obvious generalisations, in that it is not difficult to generate
sequences of weight systems of which the immanent variety form the simplest
example. For instance, one may count the number of graph morphisms from some
more sophisticated graph into the LIG, appropriately weighted. Such
generalisations would presumably filter FDs according to genus. The difficult
and interesting task is to seek topological candidates for these
generalisations. 

The Melvin-Morton-Rozansky conjecture follows naturally from our
results here (as was certainly anticipated by Bar-Natan and Garoufalidis
when they formulated their conjecture) and is the subject of a paper in
preparation. It would be interesting to more fully 
incorporate the Lie algebraic weight systems
into this picture: such an incorporation might lend some insight into
the role of Lie algebras in the space of weight systems. In \cite{BNG} an
unusual generating formula was provided for the $sl(2)$ weight system. This
formula indicates that the lesser powers in $\lambda$ arise by considering
cycle decompositions with a certain number of ``self-intersections'' in
generalised intersection graphs, with 
an extra singular point for each reduction in $\lambda$.
 
In the context of Lie algebras and cabling operators, it is appropriate to
point the reader towards \cite{AT}. In this work, Atiyah and Tall thoroughly
investigate the ``Adams operations'' on lambda rings: in our case the act of
cabling descends to such an operation on the representation ring of the 
Lie algebra. 

We finish by noting that this work relates to the BF topological field theories
investigated by Cattaneo {\it et al}. In this reference the
 authors recovered the Alexander-Conway
polynomial from the BF theory without cosmological constant \cite{CCM},
 and the Jones
polynomial from the BF theory with cosmological constant \cite{CFM}.
 Our work suggests that the 
correlators yielding the Alexander-Conway polynomial in the theory without
cosmological constant can be related to the correlators for the theory
with cosmological constant evaluated along a cabled Wilson loop. 

\
	
\noindent{\underbar{\bf Acknowledgements}}

\noindent AK was supported by an Australian Government Postgraduate Award
 and would like to express his thanks to the 
strings group at Queen 
Mary \& Westfield College London, where he was a guest whilst this work was
begun. BS is supported by the Engineering and Physical Sciences Research
Council of the UK. We would like to thank Paul Martin and Jose Figueroa-O'Farrill
for helpful conversations, and the latter also for kindly allowing us to
use his routines for drawing chord diagrams and for showing us his work prior
to publication.

\end{document}